\documentclass[sort&compress]{elsarticle}
\usepackage[normal]{subfigure}
\usepackage[margin=1in]{geometry}

\usepackage{amsmath,amssymb,amsfonts} 

\newcommand{\bs}[0]{\boldsymbol}
\newcommand{\eps}{{\displaystyle\varepsilon}}

\newcommand{\bx}[0]{\textbf{x}}
\newcommand{\bz}[0]{\textbf{z}}
\newcommand{\bp}[0]{\textbf{p}}
\newcommand{\bu}[0]{\textbf{u}}
\newcommand{\bv}[0]{\textbf{v}}
\newcommand{\by}[0]{\textbf{y}}

\newcommand{\bsub}{\begin{subequations}}
\newcommand{\esub}{\end{subequations}$\!$}

\newcommand{\revone}[1]{{\color{black}{ #1 }}}
\newcommand{\revtwo}[1]{{\color{black}{ #1 }}}
\newcommand{\revthree}[1]{{\color{black}{ #1 }}}

\newcommand{\revSec}[1]{{\color{black}{ #1 }}}

\usepackage[margin=1in]{geometry}

\usepackage{bbm}
\usepackage{graphicx}
\usepackage{pstool}
\usepackage{wrapfig}
\usepackage{caption}
\usepackage[section]{placeins} 

\usepackage{fancyhdr} 
\usepackage{lastpage} 
\usepackage{extramarks} 

\usepackage{xcolor} 
\usepackage{enumerate} 
\usepackage{paralist} 
\usepackage{amsmath, amsthm, amssymb, mathtools}
\usepackage{mathabx, pifont, stmaryrd} 
\usepackage[explicit]{titlesec} 
\usepackage{etoolbox} 
\usepackage{bibentry} 
\makeatletter\let\saved@bibitem\@bibitem\makeatother 
\usepackage[colorlinks, bookmarksopen, bookmarksnumbered,
            citecolor=red,urlcolor=red]{hyperref} 
\makeatletter\let\@bibitem\saved@bibitem\makeatother 
\usepackage{rotating}

\usepackage{algorithm}
\usepackage{algorithmic}

\theoremstyle{definition}

\usepackage{bm} 
\usepackage{amsfonts} 

\newcommand\fbm{{\ensuremath{\bm{f}}}}

\newcommand\ubm{{\ensuremath{\bm{u}}}}

\usepackage{tikz}
\usepackage{pgfplots}
\usepackage{pgfplotstable, filecontents, booktabs}
\pgfplotsset{compat=1.9}

\usetikzlibrary{pgfplots.groupplots}
\usepgfplotslibrary{fillbetween}
\usetikzlibrary{calc,fit,matrix,arrows,automata,positioning,shapes}
\usetikzlibrary{arrows.meta}

\pgfplotsset{select coords between index/.style 2 args={
    x filter/.code={
        \ifnum\coordindex<#1\fi
        \ifnum\coordindex>#2\fi
    }
}}

\tikzset{
 invisible/.style={opacity=0},
 visible on/.style={alt={#1{}{invisible}}},
 alt/.code args={<#1>#2#3}{%
   \alt<#1>{\pgfkeysalso{#2}}{\pgfkeysalso{#3}}
 },
}

\newcommand{\colorbarMatlabParula}[5]{
\begin{tikzpicture}
\begin{axis}[
   hide axis, scale only axis,
   height=0pt, width=0pt,
   colormap={parula}{rgb255=(62,38,168) rgb255=(62,39,172) rgb255=(63,40,175) rgb255=(63,41,178) rgb255=(64,42,180) rgb255=(64,43,183) rgb255=(65,44,186) rgb255=(65,45,189) rgb255=(66,46,191) rgb255=(66,47,194) rgb255=(67,48,197) rgb255=(67,49,200) rgb255=(67,50,202) rgb255=(68,51,205) rgb255=(68,52,208) rgb255=(69,53,210) rgb255=(69,55,213) rgb255=(69,56,215) rgb255=(70,57,217) rgb255=(70,58,220) rgb255=(70,59,222) rgb255=(70,61,224) rgb255=(71,62,225) rgb255=(71,63,227) rgb255=(71,65,229) rgb255=(71,66,230) rgb255=(71,68,232) rgb255=(71,69,233) rgb255=(71,70,235) rgb255=(72,72,236) rgb255=(72,73,237) rgb255=(72,75,238) rgb255=(72,76,240) rgb255=(72,78,241) rgb255=(72,79,242) rgb255=(72,80,243) rgb255=(72,82,244) rgb255=(72,83,245) rgb255=(72,84,246) rgb255=(71,86,247) rgb255=(71,87,247) rgb255=(71,89,248) rgb255=(71,90,249) rgb255=(71,91,250) rgb255=(71,93,250) rgb255=(70,94,251) rgb255=(70,96,251) rgb255=(70,97,252) rgb255=(69,98,252) rgb255=(69,100,253) rgb255=(68,101,253) rgb255=(67,103,253) rgb255=(67,104,254) rgb255=(66,106,254) rgb255=(65,107,254) rgb255=(64,109,254) rgb255=(63,110,255) rgb255=(62,112,255) rgb255=(60,113,255) rgb255=(59,115,255) rgb255=(57,116,255) rgb255=(56,118,254) rgb255=(54,119,254) rgb255=(53,121,253) rgb255=(51,122,253) rgb255=(50,124,252) rgb255=(49,125,252) rgb255=(48,127,251) rgb255=(47,128,250) rgb255=(47,130,250) rgb255=(46,131,249) rgb255=(46,132,248) rgb255=(46,134,248) rgb255=(46,135,247) rgb255=(45,136,246) rgb255=(45,138,245) rgb255=(45,139,244) rgb255=(45,140,243) rgb255=(45,142,242) rgb255=(44,143,241) rgb255=(44,144,240) rgb255=(43,145,239) rgb255=(42,147,238) rgb255=(41,148,237) rgb255=(40,149,236) rgb255=(39,151,235) rgb255=(39,152,234) rgb255=(38,153,233) rgb255=(38,154,232) rgb255=(37,155,232) rgb255=(37,156,231) rgb255=(36,158,230) rgb255=(36,159,229) rgb255=(35,160,229) rgb255=(35,161,228) rgb255=(34,162,228) rgb255=(33,163,227) rgb255=(32,165,227) rgb255=(31,166,226) rgb255=(30,167,225) rgb255=(29,168,225) rgb255=(29,169,224) rgb255=(28,170,223) rgb255=(27,171,222) rgb255=(26,172,221) rgb255=(25,173,220) rgb255=(23,174,218) rgb255=(22,175,217) rgb255=(20,176,216) rgb255=(18,177,214) rgb255=(16,178,213) rgb255=(14,179,212) rgb255=(11,179,210) rgb255=(8,180,209) rgb255=(6,181,207) rgb255=(4,182,206) rgb255=(2,183,204) rgb255=(1,183,202) rgb255=(0,184,201) rgb255=(0,185,199) rgb255=(0,186,198) rgb255=(1,186,196) rgb255=(2,187,194) rgb255=(4,187,193) rgb255=(6,188,191) rgb255=(9,189,189) rgb255=(13,189,188) rgb255=(16,190,186) rgb255=(20,190,184) rgb255=(23,191,182) rgb255=(26,192,181) rgb255=(29,192,179) rgb255=(32,193,177) rgb255=(35,193,175) rgb255=(37,194,174) rgb255=(39,194,172) rgb255=(41,195,170) rgb255=(43,195,168) rgb255=(44,196,166) rgb255=(46,196,165) rgb255=(47,197,163) rgb255=(49,197,161) rgb255=(50,198,159) rgb255=(51,199,157) rgb255=(53,199,155) rgb255=(54,200,153) rgb255=(56,200,150) rgb255=(57,201,148) rgb255=(59,201,146) rgb255=(61,202,144) rgb255=(64,202,141) rgb255=(66,202,139) rgb255=(69,203,137) rgb255=(72,203,134) rgb255=(75,203,132) rgb255=(78,204,129) rgb255=(81,204,127) rgb255=(84,204,124) rgb255=(87,204,122) rgb255=(90,204,119) rgb255=(94,205,116) rgb255=(97,205,114) rgb255=(100,205,111) rgb255=(103,205,108) rgb255=(107,205,105) rgb255=(110,205,102) rgb255=(114,205,100) rgb255=(118,204,97) rgb255=(121,204,94) rgb255=(125,204,91) rgb255=(129,204,89) rgb255=(132,204,86) rgb255=(136,203,83) rgb255=(139,203,81) rgb255=(143,203,78) rgb255=(147,202,75) rgb255=(150,202,72) rgb255=(154,201,70) rgb255=(157,201,67) rgb255=(161,200,64) rgb255=(164,200,62) rgb255=(167,199,59) rgb255=(171,199,57) rgb255=(174,198,55) rgb255=(178,198,53) rgb255=(181,197,51) rgb255=(184,196,49) rgb255=(187,196,47) rgb255=(190,195,45) rgb255=(194,195,44) rgb255=(197,194,42) rgb255=(200,193,41) rgb255=(203,193,40) rgb255=(206,192,39) rgb255=(208,191,39) rgb255=(211,191,39) rgb255=(214,190,39) rgb255=(217,190,40) rgb255=(219,189,40) rgb255=(222,188,41) rgb255=(225,188,42) rgb255=(227,188,43) rgb255=(230,187,45) rgb255=(232,187,46) rgb255=(234,186,48) rgb255=(236,186,50) rgb255=(239,186,53) rgb255=(241,186,55) rgb255=(243,186,57) rgb255=(245,186,59) rgb255=(247,186,61) rgb255=(249,186,62) rgb255=(251,187,62) rgb255=(252,188,62) rgb255=(254,189,61) rgb255=(254,190,60) rgb255=(254,192,59) rgb255=(254,193,58) rgb255=(254,194,57) rgb255=(254,196,56) rgb255=(254,197,55) rgb255=(254,199,53) rgb255=(254,200,52) rgb255=(254,202,51) rgb255=(253,203,50) rgb255=(253,205,49) rgb255=(253,206,49) rgb255=(252,208,48) rgb255=(251,210,47) rgb255=(251,211,46) rgb255=(250,213,46) rgb255=(249,214,45) rgb255=(249,216,44) rgb255=(248,217,43) rgb255=(247,219,42) rgb255=(247,221,42) rgb255=(246,222,41) rgb255=(246,224,40) rgb255=(245,225,40) rgb255=(245,227,39) rgb255=(245,229,38) rgb255=(245,230,38) rgb255=(245,232,37) rgb255=(245,233,36) rgb255=(245,235,35) rgb255=(245,236,34) rgb255=(245,238,33) rgb255=(246,239,32) rgb255=(246,241,31) rgb255=(246,242,30) rgb255=(247,244,28) rgb255=(247,245,27) rgb255=(248,247,26) rgb255=(248,248,24) rgb255=(249,249,22) rgb255=(249,251,21) },
   colorbar horizontal,
   point meta min=#1, point meta max=#5,
   colorbar style={width=10cm, xtick={#1,#2,#3,#4,#5}}
]
\addplot [draw=none] coordinates {(0,0)};
\end{axis}
\end{tikzpicture}
}

\usepackage{setspace}

\newbool{fastcompile}

\begin{document}
\title{Numerical bifurcation analysis of post-contact states in mathematical models of Micro-Electromechanical Systems}

\author[rvt1]{Charles J. Naudet\fnref{fn1}}
\ead{cnaudet@nd.edu}

\author[rvt2]{Alan E. Lindsay\fnref{fn2}\corref{cor1}}
\ead{a.lindsay@nd.edu}

\address[rvt1]{Department of Aerospace and Mechanical Engineering, University
               of Notre Dame, Notre Dame, IN 46556, United States}
\address[rvt2]{Department of Applied and Computational Mathematics and Statistics, University
               of Notre Dame, Notre Dame, IN 46556, United States}
\cortext[cor1]{Corresponding author}

\fntext[fn1]{Department of Aerospace and Mechanical
             Engineering, University of Notre Dame}
\fntext[fn2]{Department of Applied and Computational Mathematics and Statistics, University of Notre Dame}

\begin{keyword} 
Bifurcation analysis; symmetry breaking; continuation methods; singular perturbations; adaptivity; MEMS.
\end{keyword}

\begin{abstract}
This paper is a computational bifurcation analysis of a non-linear partial differential equation (PDE) characterizing equilibrium configurations in Micro electromechanical Systems (MEMS). MEMS are engineering systems that utilize electrostatic forces to actuate elastic surfaces. The potential equilibrium states of MEMS are described by solutions of a singularly perturbed elliptic nonlinear PDE. We develop a numerical method which couples a finite element approximation with mesh refinement to a pseudo arc-length continuation algorithm to numerically obtain bifurcation diagrams in the physically relevant two dimensional scenario. Several geometries, including a unit disk, square, and annulus, are studied to understand the behavior of the system over a range of domains and parameter regimes. We find that solution multiplicity, and importantly the potential for bistability in the system, depends sensitively on the parameters. In the annulus domain,  symmetry breaking bifurcations are located and asymmetric solution branches are tracked. This work significantly extends the envelope for numerical characterization of equilibrium states in microscopic electrostatic contact problems relating to MEMS. 
\end{abstract}

\maketitle

\section{Introduction}

 Micro electromechanical Systems (MEMS) are engineering systems that utilize electrostatic forces and mobile elastic components to perform a variety of tasks on miniature scales. Electrostatic actuation is effective at tiny scales due the relative importance of surface area effects over volumetric effects as linear dimension decreases. However, stable operation of MEMS are compromised when electrostatic forces overwhelm the restorative forces that elastic components can exert such that they come into contact. Predicting the onset of this criticality (known as the \lq\lq pull-in instability\rq\rq\ or \lq\lq touchdown\rq\rq) is key to MEMS design, and in particular, understanding how pull-in depends on parameter values and device geometry \cite{Pelesko2002,PB}.
 
\begin{figure}[htbp]
\centering
\subfigure[Schematic diagram.]{\includegraphics[width = 0.45\textwidth]{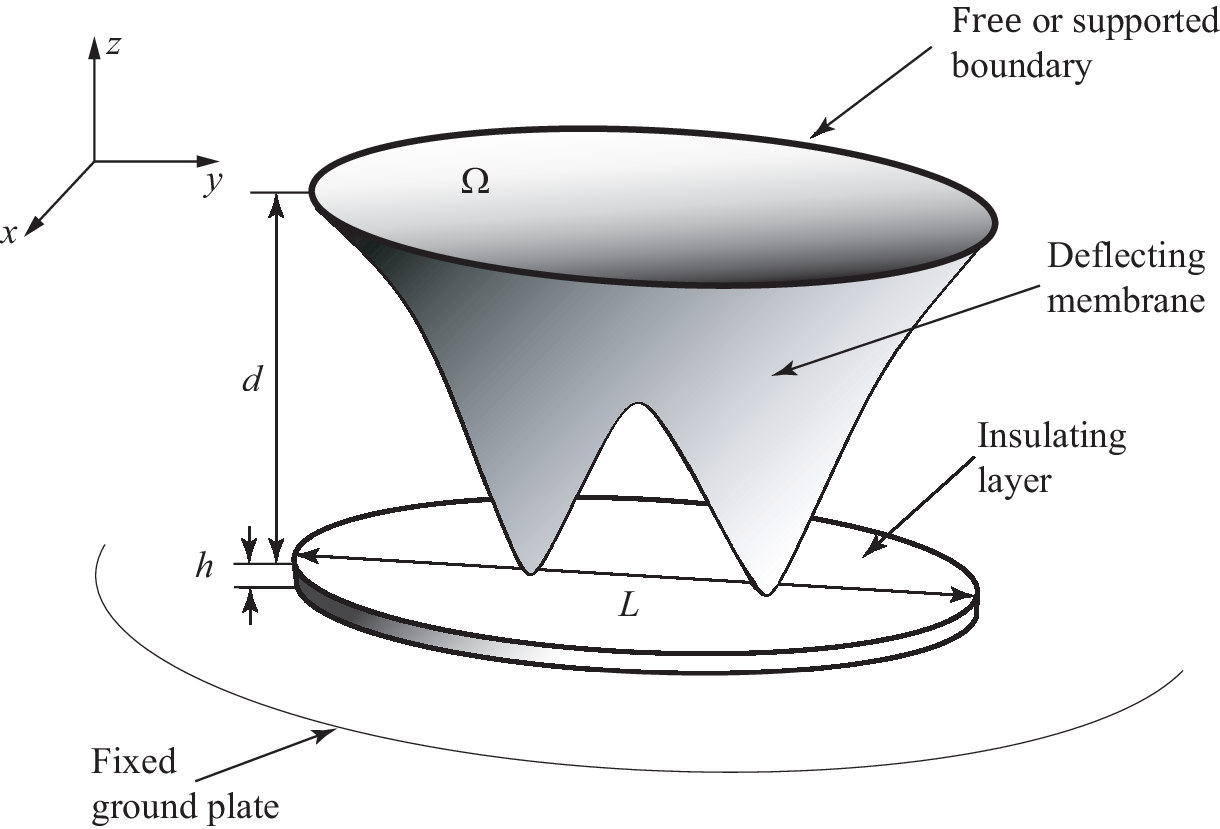} \label{fig:intro_diagram_A}}
\qquad
\subfigure[A MEMS device (source: Sandia Nat. Lab.)]{\includegraphics[width = 0.45\textwidth]{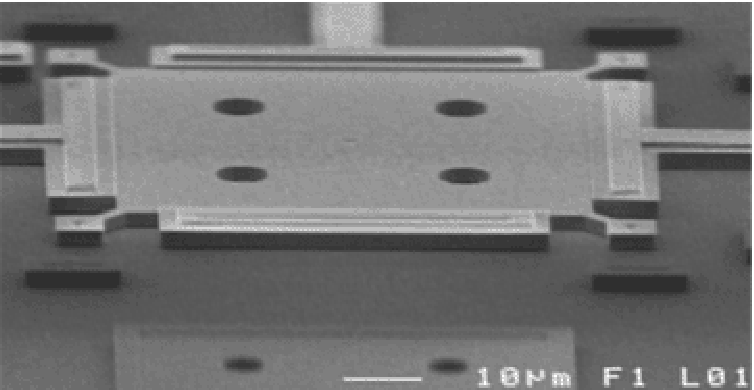} \label{fig:intro_diagram_B}}
\caption{A MEMS device (right) and a schematic (left) around which models are formulated. \label{fig:intro_diagram}}
\end{figure}

A common scenario in MEMS design features a capacitor device with a deformable elastic membrane held above an immobile substrate and fixed on its boundary. Motion of the upper surface is induced by applying an electric potential $V$ between these surfaces. If $V$ is small enough, the deflection will reach an equilibrium where the surfaces remain separated. However, if $V$ exceeds the \emph{pull-in voltage} $V^{\ast}$, the top plate will touchdown on the substrate. In Fig.~\ref{fig:intro_diagram}, a MEMS device is shown, along with a schematic representation for modeling purposes.

To describe the configuration of MEMS after touchdown, and to quantify circumstances whereby a separated state can be recovered, a family of regularized models was introduced in \cite{Equilibrium2014}. The model describes the dimensionless deflection $u(\bx,t)$ of the upper surface with the semi-linear equation 
\begin{equation}\label{mems_intro_eqns_a}
\frac{\partial u}{\partial t} = \Delta u - \frac{\lambda}{(1+u)^2} + \frac{\lambda\eps^{m-2}}{(1+u)^m},  \quad  \bx\in\Omega; \qquad u = 0, \quad \bx\in\partial\Omega,
\end{equation}
\revthree{and $\Delta = \sum_{i=1}^d \partial^2_{x_ix_i}$.}
In the formulation \eqref{mems_intro_eqns_a}, the elastic surface is a membrane occupying $\Omega$, a bounded region of $\mathbb{R}^\revtwo{d}$. The dimensions important for MEMS applications are $\revtwo{d} =1,2$. The key bifurcation parameter $\lambda\ \propto\ V^{2}$ quantifies the relative importance of electrostatic to elastic forces in the system. The small parameter $\eps$ in \eqref{mems_intro_eqns_a} controls the strength of the regularizing term $\lambda \eps^{m-2} (1+u)^{-m}$ that accounts for effects which become important when the surfaces are in close proximity. For example $m=4$ accounts for the Casimir effect (with sign of regularizing term reversed) while $m=3$ models Van der Waals forces \cite{Batra2008, GuoZhao2004}. 

The scenario $\eps=0$ is well studied \cite{PB,MEMSBook} and we briefly review some of the key properties. The most salient is the existence of a \lq\lq pull-in\rq\rq\ voltage $\lambda^{\ast}>0$ such that for all $\lambda>\lambda^{\ast}$ there are no equilibrium solutions and a finite time singularity develops such that $u(\bx,t) \to -1$ as $t\to t_c^{-}$. The limiting behavior of $u(\bx,t)$ at $t\to t_c$ has been established with estimates for both the touchdown time $t_c$ and location $\bx_c$ \cite{MEMSBook,Ward2005,LW}. For $\eps=0$ and $\lambda<\lambda^{\ast}$, the multiplicity of equilibrium solutions to \eqref{mems_intro_eqns_a} is highly dependent on domain geometry. In particular, for $2 \leq \revtwo{d} \leq 7 $, the main bifurcation branch features an infinite number of fold points implying that solution multiplicity varies considerably with $\lambda$ \cite{LW1,MEMSBook}.

\begin{figure}
\centering
\raisebox{-0.5\height}{\input{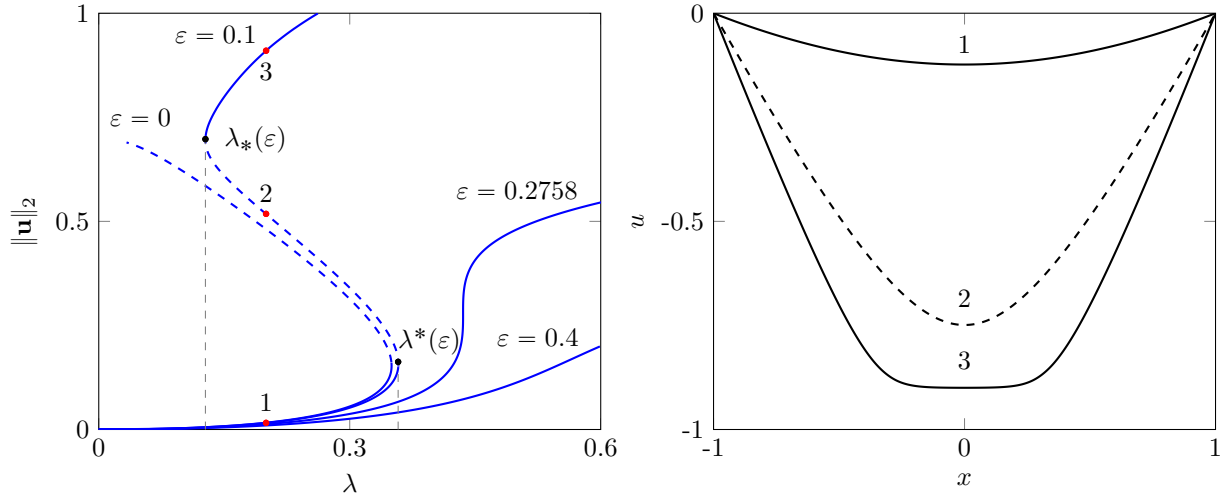}}
\caption{Numerical solutions of \eqref{eq:steady_problem} for $\Omega = [-1,1]$ and $m=4$. Left: Bifurcation curves $(\lambda,\| \bu\|_2)$ of \eqref{eq:steady_problem} for a range of $\eps$ values. For $0<\eps<0.2758$, the system is bistable in the interval $\lambda\in(\lambda_{\ast}(\eps),\lambda^{\ast}(\eps))$. \revtwo{Stable (\ref{line:stable_bif_curve_1D_eps_0.1}) branches and unstable (\ref{line:unstable_bif_curve_1D_eps_0.1}) branches are also shown.} Right: Stable \eqref{line:lam_0.2_sol1} and unstable \eqref{line:lam_0.2_sol2} solution profiles at points \eqref{line:first_sol_lamc} for $\eps = 0.1$ and $\lambda = 0.2$.\label{fig:intro1D}}
\end{figure}

In the regime $\eps>0$, touchdown is no longer a finite time singularity and solutions remain bounded \cite{Lindsay2013}. Dynamics beyond the initial contact are characterized by a sharp interface which propagates outwards and is eventually pinned by the boundary \cite{LindsayIMA2015,LD2017}. The equilibrium configurations satisfy the
equation
\begin{equation}\label{eq:steady_problem}
    \Delta u = \frac{\lambda}{(1+u)^2} - \frac{\lambda\eps^{m-2}}{(1+u)^m}, \quad \bx\in\Omega; \qquad u = 0, \qquad \bx\in\partial\Omega.
\end{equation}
%
%
The bifurcation diagram of \eqref{eq:steady_problem} for the one-dimensional domain $\Omega = [-1,1]$ reveals (cf.~Fig.~\ref{fig:intro1D}) \revthree{the presence of an imperfect bifurcation \cite{golubitsky2012singularities,ikeda2002imperfect} and a system that can exhibit bistability and hysteresis.} Specifically, for $\eps\in(0,0.2758)$ and $\lambda\in(\lambda_{\ast}(\eps),\lambda^{\ast}(\eps))$, the system has two stable solutions and one unstable. Perturbation results for the fold points $\lambda_{\ast}(\eps)$ and $\lambda^{\ast}(\eps)$ have been obtained \cite{Iuorio2019,Lindsay2016} in the limit as $\eps\to0$. Bistability allows for MEMS to switch between stable separated and contacted states and hence can be utilized for practical effect \cite{Batra2008,bistability2011}. Consequently, there has been significant interest in modeling these contact states together with analysis of their dynamics and stability \cite{GOLDBERG2022,Walker2021}. Of particular practical importance is determining the extent of bistable operation which requires determining the equilibrium structure of \eqref{mems_intro_eqns_a} as the physical parameters $\lambda$, $\eps$ vary \cite{Lindsay2016,Iuorio2019}. In higher dimensional settings, even in radially symmetric cases, little is currently known about the bifurcation structure of \eqref{eq:steady_problem}. The system \eqref{eq:steady_problem} belongs to a class of challenging singular perturbation problems featuring two independent small parameters \cite{KUEHN2022133105}.

In the present work, we are focused on numerical calculation of the bifurcation diagram of \eqref{eq:steady_problem} in two dimensions, particularly in the singularly perturbed and physically relevant case $\eps\ll1$. In this regime, we find that equation \eqref{eq:steady_problem} has a remarkably complex solution structure. For \eqref{eq:steady_problem} in the 1D domain $\Omega = [-1,1]$, both the numerical treatment and the observed solution multiplicity are significantly simpler. Computationally, there are many \lq\lq tricks\rq\rq\ and reformulations which do not extend naturally to higher dimensions. For example, to calculate the 1D bifurcation diagram of \eqref{eq:steady_problem} shown in Fig,~\ref{fig:intro1D}, we define for $\by = (y_1,y_2,y_3,y_4) = (u,u_x,\lambda,\int_{-1}^x u^2dx)$, the augmented boundary value problem (BVP) 
\begin{equation}
\frac{d}{dx}\begin{pmatrix}
    y_1\\y_2\\y_3\\y_4
\end{pmatrix} = \begin{pmatrix}
    y_2\\ y_3\big[ \frac{1}{(1+y_1)^2} -\frac{\eps^2}{(1+y_1)^{m}}\big] \\ 0 \\ y_1^2
\end{pmatrix}, 
\qquad
\begin{array}{r}
y_1(-1) = 0\\
y_1(1) = 0\\
y_4(-1) = 0\\
y_4(1) = \sigma
\end{array}.
\end{equation}
This system is then solved for prescribed values of $\sigma=\|u\|^2_2$ with a standard BVP routine (e.g. {\tt bvp4c} in \textsc{Matlab}). This approach assumes the system can be reformulated into an ODE and that each solution is uniquely prescribed by its $\|u\|_2$ value. In higher dimensions, neither of these conditions are necessarily met.

The contribution of the current work is to combine adaptive finite element simulations \cite{funken2011efficient} with pseudo arc-length continuation \cite{Allgower03,Krauskopf07} to numerically compute bifurcation diagrams for \eqref{eq:steady_problem} in general two dimensional geometries. There are two main challenges that arise. First, the multiplicity of solutions depends sensitively on the parameters $(\lambda,\eps)$. To address this, we apply a finite element method (FEM) discretization of \eqref{eq:steady_problem} to obtain a nonlinear system of general form
\begin{equation}\label{eqn:discrete}
    \fbm(\bu,\bs{\alpha}) = 0, \qquad \fbm: \mathbb{R}^{\revtwo{N_\ubm}} \times \mathbb{R}^z \to \mathbb{R}^\revtwo{N_\ubm}.
\end{equation}
Here $\bs{\alpha}\in\mathbb{R}^z$ represents the vector of parameters \revtwo{and $N_\ubm$ represents the number of degrees of freedom to represent the solution.} We are principally interested in characterizing the solution set over the parameter set $\bs{\alpha} = (\lambda,\eps)$. In Sec.~\ref{sec:newt_raph} we describe an implementation of the pseudo-arc length continuation method \cite{Allgower03,arclength2001,TataKeller,Chan1982} to trace out solution branches of \eqref{eqn:discrete} with sharp snap backs \cite{hellweg1998new}. The second challenge is that the equilibrium solutions to \eqref{eq:steady_problem}, particularly those where $u\approx -1$, have localized sharp features. To address this, we employ adaptive mesh refinement \cite{funken2011efficient} within our arc-length continuation method to resolve \revthree{rapidly changing solution features,}thereby allowing for accurate continuation along singular solution branches. This work therefore represents a significant extension upon previous numerical methodologies for bifurcation analysis and allows for exploration of solution structures that were previously only available in radially symmetric scenarios and through exploitation of symmetries \cite{PB} or reformulation into initial value problems \cite{LW1}. \revtwo{The numerical experiments in this work explore the case $m=4$ and we observe that parameter values $m>2$ yield qualitatively similar behaviors.}

The structure of the paper is the following. In Sec.~\ref{sec:numMethods} we outline our numerical method which combines an adaptive finite element method with pseudo arc-length continuation. In Sec.~\ref{sec:Results} we show results of the method after initial validation on a previously solved case (Sec.~\ref{sec:conv_study}). Following this, we solve the main equation \eqref{eq:steady_problem} for values $\eps = [0.3,0.1,0.01]$ on unit disk and square geometries where we observe highly sensitive dependence of solution multiplicity on parameters $(\eps,\lambda)$. Generically we observe that as $\eps$ decreases, the solution structure becomes more intricate with multiplicity depending sensitively on $\lambda$. Following on, we investigate in Sec.~\ref{sec:symm_break} the non-simply connected scenario in the case of an annular domain $r_1 < |\bx|<1$. In the scenario $\eps=0$, we numerically calculate the first nine of a sequence of symmetry breaking bifurcations on the radially symmetric solution branch, each of which appears to give rise to an infinite fold point structure. The resulting bifurcation diagram (cf.~Fig.~\ref{fig:annulus_bif_curves}) shows an extremely intricate solution structure, with particularly sensitive solution multiplicity as $\lambda\to0$. \revthree{When $\eps>0$ (cf.~Fig.~\ref{fig:asymBranches}), we find that the radially symmetric branch now has a finite number of symmetry breaking bifurcations.} Finally, in Sec. \ref{sec:discussion} we discuss some avenues for future investigations.

\section{Numerical methods for discretization and arc length continuation}\label{sec:numMethods}
\subsection{Continuous Galerkin discretization of steady model} \label{sec:fem}

Our numerical solution of \eqref{eq:steady_problem} is based on a finite element discretization \cite{funken2011efficient} on a triangulation $\mathcal{T}_h$ of a two-dimensional domain $\Omega_h$ that approximates $\Omega$. The spatial resolution of the triangulation can be adapted to resolve sharp features ($h$-adaptive). We adopt a linear finite \revone{trial and test} space $V_h$ which is continuous in $\Omega_h$ and vanishes on $\partial\Omega_h$. Specifically,
\begin{equation} \label{eq:linear_poly_space}
    V_h = \left\{ u_h \in C(\Omega_h)  \ \mid \ u_h \mid_{\partial\Omega_h} = 0, \quad u_h\mid_{T} \ \in \mathbb{P}_1, \quad \forall T \in \mathcal{T}_h \ \right\},
\end{equation}
where $\mathbb{P}_1$ is the set of linear polynomials . From the vertices $\{\bp_j\}_{j=1}^{\revtwo{N_\bp}}$ of $\mathcal{T}_h$, the basis $\phi_i\in\mathcal{V}_h$ is constructed so that $\phi_i(\bp_j) = \delta_{ij}$.\revone{Multiplying \eqref{eq:steady_problem} by the test function and integrating by parts, we obtain the weak form over a given element $T \in \mathcal{T}_h$: find $u_h \in V_h$ for all $\phi_h \in V_h$
\begin{equation}
    \int_{T}\nabla \phi_h \cdot \nabla u_h\, dV +  \int_{T}  \phi_h \cdot \mathcal{S}(u_h, \lambda)\, dV = 0,
\end{equation}
where $\mathcal{S}(u_h, \lambda)$ is the source term on the right hand side of \eqref{eq:steady_problem} defined as 
\begin{equation}
    \mathcal{S}(w, \lambda) = \frac{\lambda}{(1+w)^2} - \frac{\lambda\eps^{m-2}}{(1+w)^m}.
\end{equation}
The residual over an element is given as 
\begin{equation}
    f_{T} : (\phi_h, u_h; \lambda) \mapsto \int_{T}\nabla \phi_h \cdot \nabla u_h\, dV +  \int_{T}  \phi_h \cdot \mathcal{S}(u_h, \lambda)\, dV = 0,
\end{equation}
and the global residual is defined in terms of element residuals as 
\begin{equation}
    f_h : (\phi_h, u_h; \lambda) \mapsto \sum_{T \in \mathcal{T}_h} f_{T}(\phi_h, u_h; \lambda).
\end{equation}
Finally, after introducing the basis associated with trial and test spaces \eqref{eq:linear_poly_space}, assumed to be the same in this work, we obtain a system of nonlinear algebraic equations
\begin{equation}
    \fbm : \mathbb{R}^{N_\ubm} \times \mathbb{R} \rightarrow \mathbb{R}^{N_\ubm} \qquad\fbm : (\ubm, \lambda) \mapsto \fbm(\ubm, \lambda),
\end{equation}
where $N_\ubm = \dim(V_h)$, which is the standard Galerkin finite element residual.
}


\subsection{Adaptive mesh refinement} \label{sec:adapt_href}
We implement an adaptive mesh refinement algorithm \cite{funken2011efficient} that refines and coarsens the mesh at each steady state, or at every $\nu$ steady state solutions, and if the solution struggles to converge, within each steady state solve. We use newest vertex bisection\revthree{where if an element is tagged for refinement, its longest edge is bisected by inserting a new node. A new edge connecting the newest vertex and the vertex opposite the longest edge is inserted, resulting in two new elements from the single parent element. This procedure alone can result in a hanging node, so the procedure is repeated on the element opposite the longest edge. Other types of refinement may be done if a given element has two or even three of its edges marked for refinement. In this case, newest vertex bisection is performed first on the longest edge while the other two edges become the reference edges in the resulting elements. Newest vertex bisection is then performed on these edges if they were marked for refinement, as well on any resulting hanging nodes. The coarsening strategy coarsens the mesh by removing nodes inserted using newest vertex bisection and merging two brother elements that resulted from newest vertex bisection into one element.} 

Elements are tagged for refinement or coarsening based on a hierarchical error estimator \cite{funken2011efficient} from the computed solution. In this way, elements are added or removed based on the error they contain. At each step of the method we compute an error vector $\textbf{e} = (e_1, \ldots,e_{N_T})$ where $N_T$ is the number of elements and define tunable parameters for refinement, $\kappa>0$, and coarsening, $\rho>0$, to modulate the refinement and coarsening strategies respectively. For each element $i$ such that \revtwo{$e_i > \kappa \sum_{j=1}^{N_T} e_j$, the element is selected for refinement. Conversely, if $e_i$ satisfies $e_i < \rho \sum_{j=1}^{N_T} e_j$} then the element is selected for coarsening. These decisions allow us to determine where the error is distributed in the domain and attribute more or less elements to regions of the mesh accordingly. \revSec{Finally, after the mesh has been coarsened or refined, the solution is interpolated from the old mesh onto the new mesh, in order to implement the pseudo-arc-length continuation technique using the same collocation points, as required.}

The algorithm increases and decreases local mesh density in regions where the solution has large and small spatial variations, respectively. Localized mesh refinement becomes essential in the regime $\eps \ll 1$ where solutions along the upper branch of the bifurcation curve develop sharp boundary layers.

\subsection{Pseudo arc-length continuation}\label{sec:newt_raph}

Obtaining the full parametric solution multiplicity of \eqref{eqn:discrete} is a challenging general problem because for any particular value of $\lambda$, there can be any number of solutions. We focus on a fixed value of $\eps$ and reduce to the scalar parameter case $\fbm(\bu,\lambda) = 0$. We now describe an implementation of pseudo arc-length continuation \cite{BEYN2002149,TataKeller} that can reveal the structure of equilibrium solutions. Starting from an initial solution $(\bu_0,\lambda_0)$ and tangent vector $(\dot{\bu}_0,\dot{\lambda}_0)$, the method seeks to determine a sequence of points $(\bu_j,\lambda_j)$ that satisfies $\fbm(\bu, \lambda) = 0$ and can traverse simple fold points. In terms of an arc-length parameter $s$, the solution curve $(\bu(s),\lambda(s))$ is augmented by solving the equation
\begin{equation}
    n(\bu, \lambda, ds) = \dot{\bu}_j ( \bu - \bu_j) + \dot{\lambda}_j (\lambda - \lambda_j) - ds,
\end{equation}
which describes the equation of a hyperplane with normal vector $(\dot{\bu}_j,\dot{\lambda}_j)$ and perpendicular distance $ds$ from $(\bu_j,\lambda_j)$. The bifurcation curve is traversed by solving the augmented system
\bsub\label{eq:arc_length}
\begin{align}
    \label{eq:arc_length_a} \fbm(\bu, \lambda) =  0,\\[4pt]
    \label{eq:arc_length_b} n(\bu,\lambda,ds) = 0,
\end{align}
\esub
with the arc-length parameter $ds$ chosen in an adaptive manner to adequately resolve the solution curve while maintaining efficiency. Specifically, we apply a Newton-Raphson method to \eqref{eq:arc_length} with the initial condition $(\bu_j,\lambda_j)$ which yields the iteration scheme

\begin{equation}\label{eq:NR_system}
\begin{pmatrix} 
\fbm_u (\bu^{(k)}, \lambda^{(k)}) & \fbm_{\lambda} (\bu^{(k)}, \lambda^{(k)})  \\[5pt]
\dot{\bu}_j & \dot{\lambda}_j
\end{pmatrix} 
\begin{pmatrix}
\Delta \bu\\[5pt] \Delta \lambda 
\end{pmatrix}
= -
\begin{pmatrix}
\fbm(\bu^{(k)}, \lambda^{(k)})\\[5pt]
n(\bu^{(k)}, \lambda^{(k)},ds)
\end{pmatrix}
\hspace{0.5in}
\begin{matrix}
\bu^{(k+1)} = \bu^{(k)} + \Delta \bu\\[5pt]
\lambda^{(k+1)} = \lambda^{(k)} + \Delta \lambda
\end{matrix}
\end{equation} 
Provided $ds$ is not too large, the iterative scheme \eqref{eq:NR_system} convergences rapidly. The tangent vector $(\dot{\bu},\dot{\lambda})$ is found by differentiation of $\fbm= 0$ with respect to $s$ to obtain the linear system $\fbm_u \bz + \fbm_{\lambda}=0$. To obtain the tangent vector, we apply the scaling
\begin{equation}
    \dot{\bu} = a\, \bz, \qquad \dot{\lambda} = a, \qquad a = \frac{\pm1}{\sqrt{1+\|\bz\|^2} }.
\end{equation}
where the sign of $a$ is chosen to preserve direction of travel along the bifurcation curve. The stepsize $ds$ is varied continuously based on the convergence of the iterative scheme \eqref{eq:NR_system} and progress along the bifurcation curve. At each admitted point on the bifurcation curve, we calculate the error estimates (Sec.~\ref{sec:newt_raph}) and refine or coarsen where necessary.

\subsection{Branch point detection and branch switching} \label{sec:branch_switch}
For families of solutions that contain symmetry breaking solutions it can be beneficial to locate the point at which there is a bifurcation and a split into different branches, then travel along those branches \cite{hellweg1998new,shi1992simple}. To detect branch points of our system, we monitor the determinant of the Jacobian, det$(\fbm_u(\bu^{(k)}, \lambda^{(k)}))$ and \revSec{observe for zeros}. A sign change in this quantity signifies a bifurcation, however, a double root whereby two eigenvalues cross zero, results in a determinant only touching zero \cite{golubitsky2012singularities,ikeda2002imperfect}. To precisely locate branch points, we stop path tracking when the determinant decreases below a pre-specified tolerance and solve an extended system \cite{Allgower03,BEYN2002149} for the non-trivial null space of the Jacobian. We introduce $\bv \in \mathbb{R}^{N_\ubm}$ and $\beta \in \mathbb{R}$ and solve 
\begin{equation}\label{eq:branch_point_locator}
  \begin{cases}
    \fbm(\bu, \lambda) + \beta \bv = 0       \\
    \fbm^T_u (\bu, \lambda) \bv = 0         \\
    \bv^T \fbm_{\lambda}(\bu,  \lambda) = 0  \\
    \bv^T \bv - 1 = 0           
  \end{cases}
\end{equation} 

This problem is solved using the Newton-Raphson method with initial condition $(\bu, \lambda, \beta, \bv) = (\bu^{(k)}, \lambda^{(k)}, 0, \bar{\bv})$ and $\bar{\bv}$ is initially defined as the eigenvector of $\fbm^T_\bu$ corresponding to the smallest real eigenvalue $\mu$, satisfying $\fbm^T_\bu \bar{\bv} = \mu \bar{\bv}$. A simple branch point $(\bu, \lambda)$ corresponds to a regular solution $(\bu, \lambda, 0, \bv)$ of \eqref{eq:branch_point_locator}.

\section{Results}\label{sec:Results}

In this section, the numerical methods presented are first tested against a known study, and then several numerical experiments are conducted that demonstrate the methods' robustness in exploring the stability and multiplicity of various problems, specifically \eqref{eq:steady_problem} that models equilibrium configurations for $\eps \geq 0$. This model will be explored for various domains $\Omega$, including a circular disk, a square, and an annulus. 

\subsection{Convergence study with circular domain example} \label{sec:conv_study}

\begin{figure}[t!]
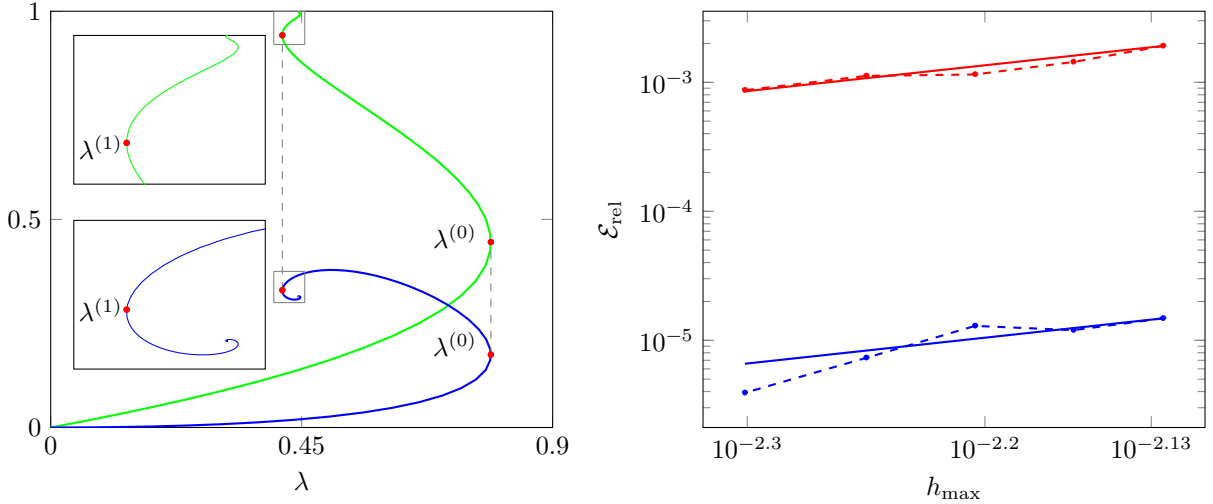

    \centering

    \raisebox{-0.5\height}{\input{_py/spiral_plot.tikz}}
    \begin{picture}(0,0)
        \put(-435, 25){\input{_py/spiral_plot_insetz.tikz}}
        \put(-435, -45){\input{_py/u2_spiral_plot_insetz.tikz}}
    \end{picture}
    
    \caption{Radially symmetric solution of \eqref{eq:steady_problem} for $\eps= 0$. Left: Bifurcation curves of $\lambda$ against $\|\bu\|_2$ \eqref{line:u2_spiral_curve} and $\|\bu\|_\infty$ \eqref{line:spiral_curve} with insets showing spiral. \revtwo{The points for respective $\lambda^{(k)}$ values are also shown on each curve \eqref{line:lam1}.} Right: Convergence of error in computation of the first two fold points $\lambda^{(0)}$ (\ref{line:e1_conv})  and $\lambda^{(1)}$ (\ref{line:e2_conv}) for a range of fixed mesh sizes $h_{\text{max}}$. We observe the expected second order convergence (\ref{line:line1_conv} and \ref{line:line2_conv}) (right).}
    \label{fig:MEMSDisk}
\end{figure}

As validation of our numerical method, we first consider the case $\eps = 0$ which has been the topic of many studies \cite{Pelesko2002,LW}. In the case of the unit disk domain $\Omega = \{ \bx \in \mathbb{R}^2 \, \mid \, |\bx| \leq 1 \}$, and applying radial symmetry $u = u(r)$ for $r = |\bx|$, equation \eqref{eq:steady_problem} reduces to the ODE
\begin{equation}
    \frac{\partial^2 u}{\partial r^2} +  \frac{1}{r} \frac{\partial u}{\partial r} = \frac{\lambda}{(1+u)^2}, \quad 0<r<1; \qquad u'(0) = 0, \quad u(1) = 0.
\end{equation}
The solutions of this boundary value problem can be mapped to an initial value problem by exploiting scale invariance \cite{Pelesko2002}. Specifically, by introducing the variables $t$ and $w$ such that
\[
u(r) = -1 + \alpha w(\eta), \qquad \eta = t r,
\]
some simple algebra shows that the bifurcation diagram can be parameterized as
\begin{equation}\label{eq:scale_bif}
|u(0)| = 1 - \frac{1}{w(\eta)}, \qquad \lambda = \frac{\eta^2}{w^3(\eta)},
\end{equation}
where $w(\eta)$ solves the parameter-free initial value problem
\begin{equation}\label{eq:scale_inv}
    w'' + \frac{1}{\eta} w' = \frac{1}{w^2}, \quad \eta >0; \qquad w(0) = 1, \quad w'(0) = 0.
\end{equation}
Therefore by solving \eqref{eq:scale_inv} to high accuracy, the bifurcation diagram is reproduced from \eqref{eq:scale_bif} as shown in Fig.~\ref{fig:MEMSDisk}. A notable feature is the presence of an infinite fold point structure which highlights the sensitive solution multiplicity as $\lambda$ varies \cite{LW1}. As validation for our arc-length continuation method, we use the first two folds $\lambda^{(0)}$ and $\lambda^{(1)}$ as a comparison. These fold points are found to high accuracy by applying a root solving routine to $\lambda_\eta = \frac{2\eta}{w^3}[1 - \frac{3}{2}\eta \frac{w'}{w}]= 0$. We obtain that
\[
\lambda^{(0)} = 0.78922927, \qquad \lambda^{(1)} = 0.41533025.
\]
To confirm the validity of our numerical method, we use a range of fixed meshes and calculate the error term
\begin{equation}\label{eq:error}
\mathcal{E}_{\textrm{rel}} = \frac{|\lambda_h^{(k)} - \lambda^{(k)}|}{\lambda^{(k)}},
\end{equation}
to observe convergence. In Fig.~\ref{fig:MEMSDisk} we plot the bifurcation curves of $\lambda$ against both $\|u\|_2$ and $\|u\|_{\infty}$ together with convergence of the relative errors \eqref{eq:error}. The characteristic feature of the bifurcation is a tight spiral with $\lambda\to\lambda_c$ as $\|u\|_{\infty}= |u(0)| \to 1$. This limiting value can be established through a re-scaling \cite{Pelesko2002} of equation \eqref{eq:scale_inv}. Setting $\xi = \log \eta$ and $w = \eta^{\frac23} v(\xi)$ yields
\begin{equation}\label{eq:rescaling}
    \frac{d^2v}{d\xi^2}+ \frac{1}{3}\frac{dv}{d\xi}+ \frac{2}{9}v = \frac{1}{v^2}.
\end{equation}
A reformulation of equation \eqref{eq:rescaling} into a planar dynamical system in variables $(1/v,v_{\xi}/v)$ yields the behavior of $v(\xi)$ as $\xi\to\infty$ which is characterized by a stable spiral which centers on $\lambda\to\lambda_c=4/9$ as $|u(0)|\to 1$. This analysis for the disk case is useful for validating our finite element approach, however, it does not generalize to non-radially symmetric domains where a fully numerical approach is needed.
\begin{figure}[t!]
\centering
 \resizebox{0.85\textwidth}{!}{\input{_py/circle_all_bif_curves.tikz}}
 \resizebox{0.85\textwidth}{!}{\begin{tikzpicture}
\begin{groupplot} [
group style={group size = 3 by 2, horizontal sep = 0.05cm, vertical sep = 0.2cm},
title style={at={(current bounding box.north)}, anchor=south}]
\nextgroupplot[axis equal image, width=0.45\textwidth, xtick={-0, 1}, ytick={}, title={1}, xticklabels={}, yticklabels={}, ylabel={$\varepsilon = 0.3$}, xlabel={}, xmin=0, xmax=1, ymin=-1, ymax=0]
\addplot []
graphics [xmin=0,xmax=1,ymin=-1,ymax=0] { 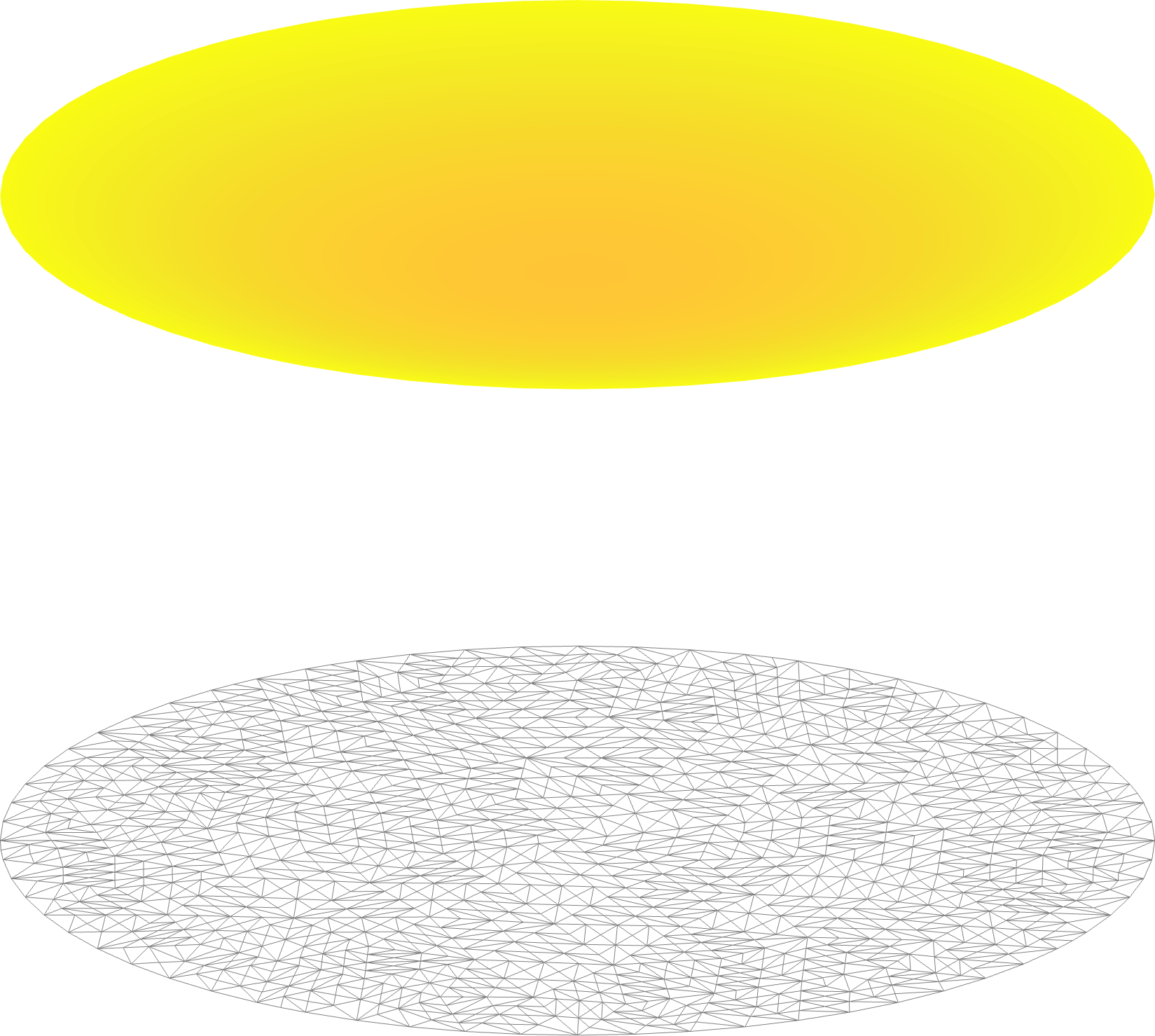};

\nextgroupplot[axis equal image, width=0.45\textwidth, xtick={}, ytick={7}, title={2}, xlabel={}, yticklabels={}, xticklabels={}, xmin=0, xmax=1, ymin=-1, ymax=0]
\addplot []
graphics [xmin=0,xmax=1,ymin=-1,ymax=0] { 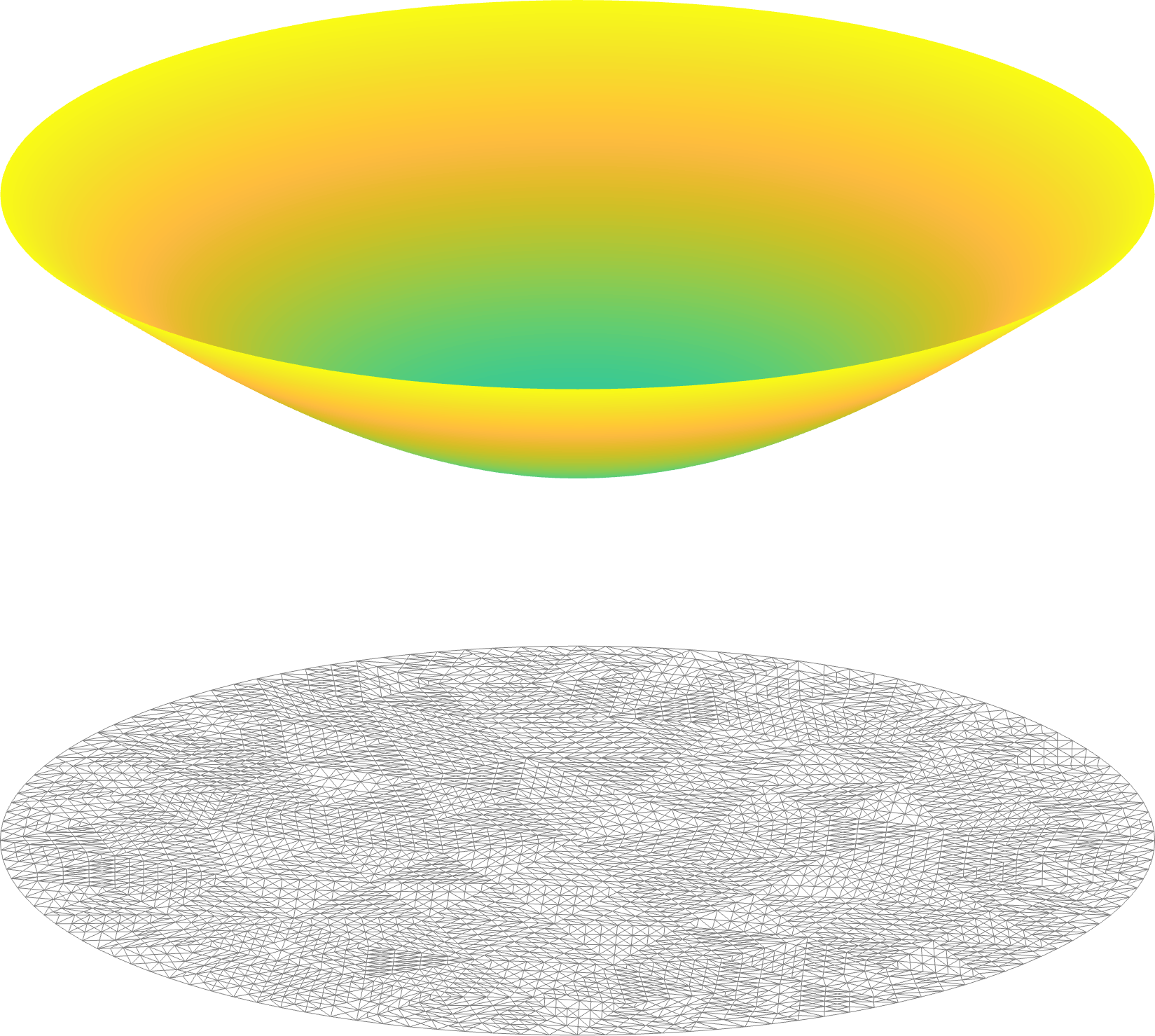};

\nextgroupplot[axis equal image, width=0.45\textwidth, xtick={}, ytick={}, title={3}, xlabel={}, xticklabels={}, yticklabels={}, xmin=0, xmax=1, ymin=-1, ymax=0]
\addplot []
graphics [xmin=0,xmax=1,ymin=-1,ymax=0] { 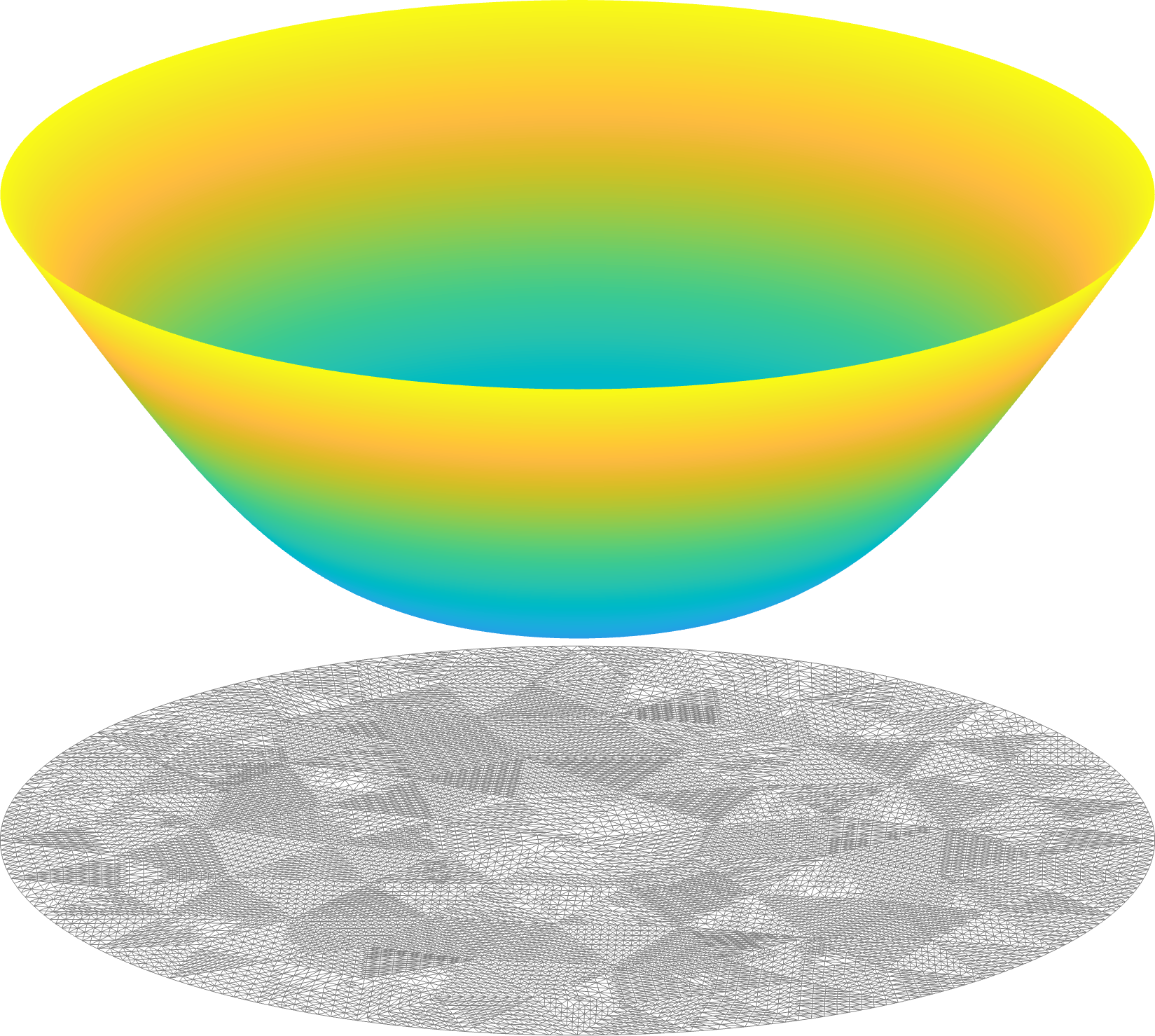};

\nextgroupplot[axis equal image, width=0.45\textwidth, xtick={}, ytick={}, xticklabels={}, yticklabels={}, ylabel={$\varepsilon = 0.1$}, xlabel={}, xmin=0, xmax=1, ymin=-1, ymax=0]
\addplot []
graphics [xmin=0,xmax=1,ymin=-1,ymax=0] { 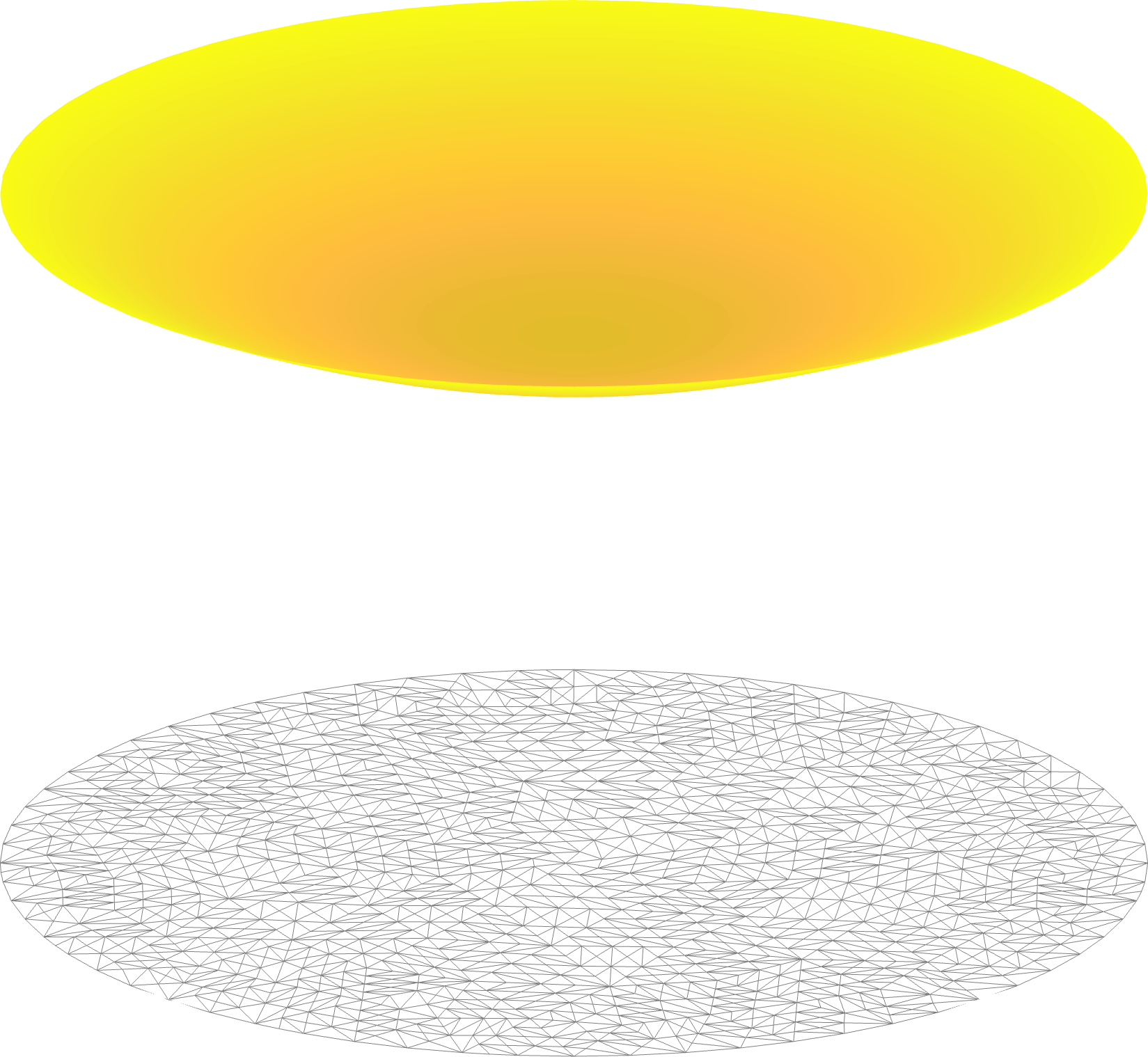};

\nextgroupplot[axis equal image, width=0.45\textwidth, xtick={}, ytick={}, xlabel={}, yticklabels={}, xticklabels={}, xmin=0, xmax=1, ymin=-1, ymax=0]
\addplot []
graphics [xmin=0,xmax=1,ymin=-1,ymax=0] { 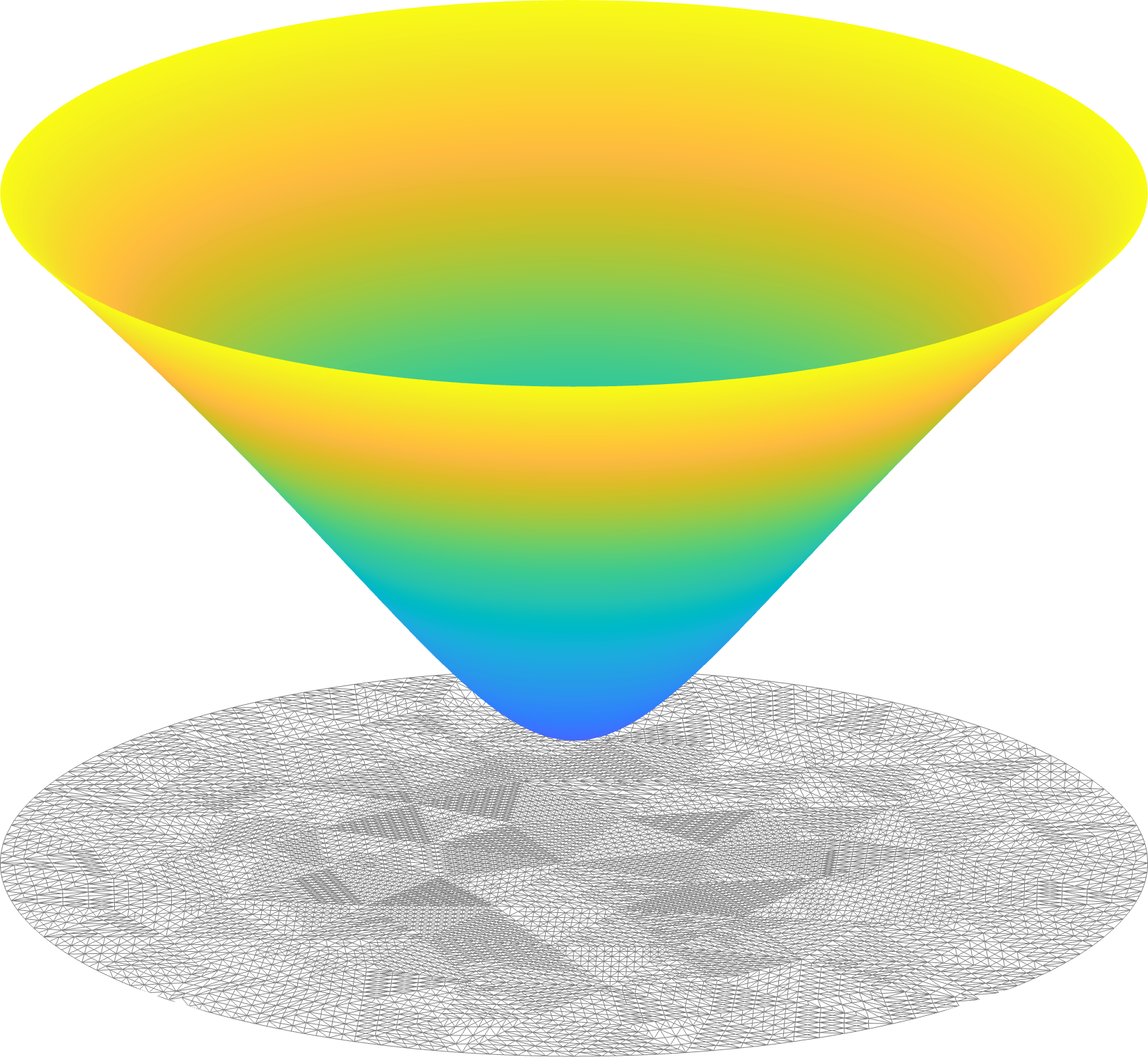};

\nextgroupplot[axis equal image, width=0.45\textwidth, xtick={}, ytick={}, xlabel={}, xticklabels={}, yticklabels={}, xmin=0, xmax=1, ymin=-1, ymax=0]
\addplot []
graphics [xmin=0,xmax=1,ymin=-1,ymax=0] { 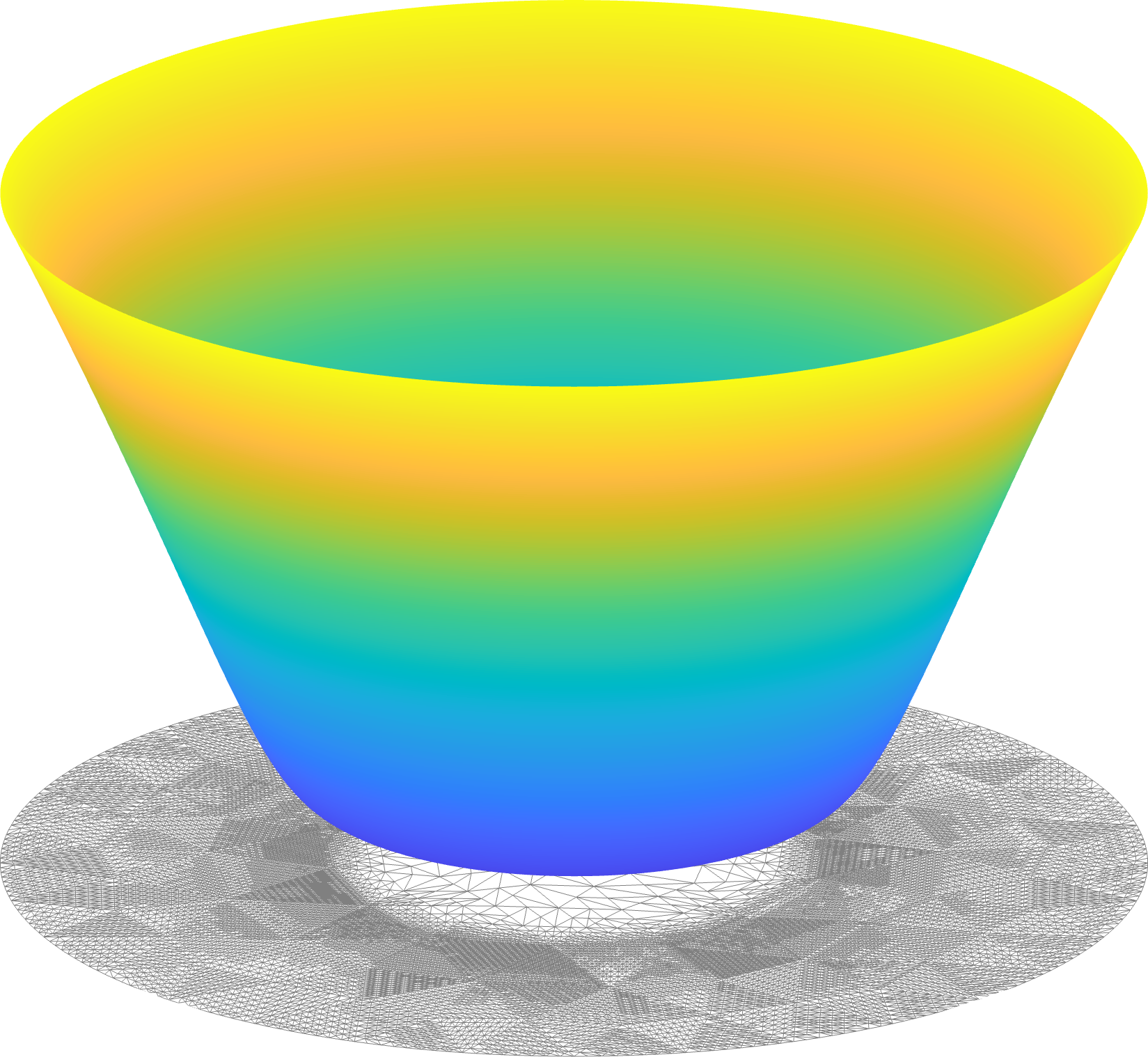};

\end{groupplot}\end{tikzpicture}}
 \colorbarMatlabParula{-1}{-0.75}{-0.5}{-0.25}{0}
 \caption{Solutions to \eqref{eq:steady_problem} on the unit disk with $m=4$. Top Row: Bifurcation curves of $\lambda$ against $\|\bu\|_2$ (\ref{line:bif_curve_eps_0.0225}) and $\|\bu\|_{\infty}$ (\ref{line:bif_curve_eps_0.1}) for various $\eps$ values. We identify two cubic fold points (\ref{line:point_eps_0.25_lam_crit}) through which solution multiplicity increases. Middle row: Solution profiles $u$ and adapted mesh for the points 1,2,3 (\ref{line:admitted_sol_eps_0_3}) on the $\eps = 0.3$ curve and $\lambda = [0.5, 1, 1.5]$. Bottom row: Solution profiles $u$ and adapted mesh for points 1,2,3 (\ref{line:admitted_sol_eps_0_3}) on the $\eps = 0.1$ curve and $\lambda = 0.6$.}
 \label{fig:circle_side_eps_0_1_lam_0_6}
\end{figure}

\subsection{Regularized problem on unit disk} \label{sec:reg_unit_disk}

In this section, results for the regularized problem \eqref{eq:steady_problem} on the unit disk $\Omega = \{ \bx \in \mathbb{R}^2 \, \mid \, |\bx| \leq 1 \}$ are presented for a range of $\eps$ values and $m=4$. We observe in Fig.~\ref{fig:circle_side_eps_0_1_lam_0_6} the bifurcation curves of $\lambda$ against $\| \bu \|_{\infty}$ and $\| \bu \|_{2}$ which show the multiplicity of solutions changing considerably as $\eps\to0$. The curve $(\lambda,\| \bu \|_{2})$ is more useful for differentiating between solutions which have similar values of $\|\bu\|_{\infty}$. We will observe in several cases the effectiveness of the adaptive refinement implementation (Sec.~\ref{sec:adapt_href}) and the pseudo arc-length continuation strategy (Sec.~\ref{sec:newt_raph}), particularly when $\eps\ll1$. We remark upon the changes in solution structure as $\eps$ decreases in value.

\underline{Case $\eps = 0.3$:} At this relatively large value, there exists a unique solution for each value of $\lambda$. The bifurcation curve is traced out from the initial point $(\lambda,\bu) = (0,\textbf{0})$ which is the first steady state in the family of solutions for \eqref{eq:steady_problem}. In Fig.~\ref{fig:circle_side_eps_0_1_lam_0_6} we show bifurcation curves (top row) of $\lambda$ against $\|\bu\|_2$ and $\|\bu\|_{\infty}$ together with three representative deflection profiles and adapted meshes (bottom two rows). We observe the deflection $u$ increasing as $\lambda$ increases with the largest point of deflection occurring at the center of the domain. Once the center of the disk attains the maximum deflection $\|\bu\|_{\infty} \approx 1- \eps = 0.7$, the source terms $\lambda(1+u)^{-2}$ and $\lambda \eps^2(1+u)^{-4}$ are balanced and further increasing $\lambda$ results in the deflection spreading radially outwards and an expansion of the flat central region where $\min_{\bx\in\Omega}u(\bx) \approx -1 + \eps = -0.7$. 

\begin{figure}[h]
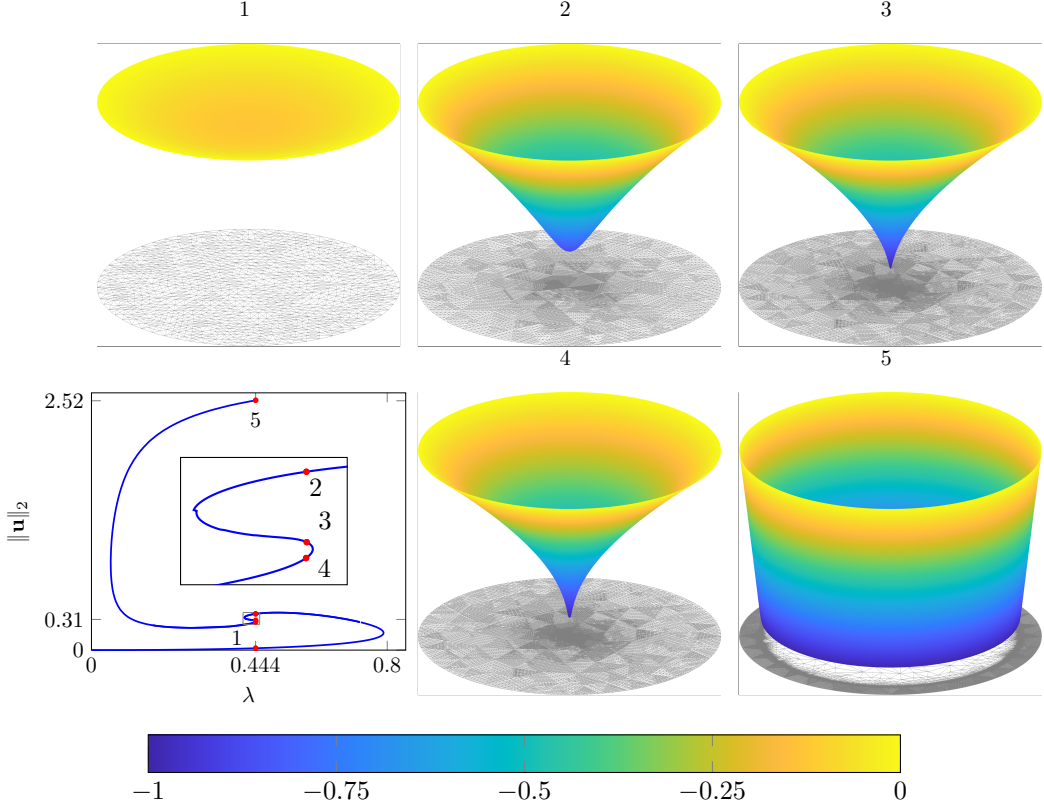

\centering
  \resizebox{0.85\textwidth}{!}{\input{_py/circle_tilted_w_msh_eps_0_01_lam_0_444.tikz}}
 \colorbarMatlabParula{-1}{-0.75}{-0.5}{-0.25}{0}
  \begin{picture}(0,0)
 \put(-285, 85){\input{_py/circle_tilted_eps_0_01_lam_0_444_inset.tikz}}
\end{picture}
 \caption{Solutions to \eqref{eq:steady_problem} on the unit disk for $\eps = 0.01$, $m=4$ with deflection profile $u$ and adapted meshes. Five solutions (panels 1-5) are shown for $\lambda = 0.444$ with their positions (\ref{line:point_eps_0.01_lam_0.444_1}) along the bifurcation curve (\ref{line:bif_curve_eps_0.01_lam_0.444}).}
 \label{fig:circle_tilted_eps_0_01_lam_0_444}
\end{figure}

As $\eps\to0$, we expect solutions with $\|\bu\|_{\infty} \approx 1-\eps$ for which the forcing terms become $\mathcal{O}(\eps^{-2})$. Correspondingly, we must resolve solutions with sharp solution features at localized regions of the domain. As $\eps$ decreases, we identify a critical bifurcation value $\eps_{c_1} = 0.25966$ for which the solution multiplicity changes from one to three. For $\eps > \eps_{c_1}$, there is a unique solution for each $\lambda$ while for $\eps < \eps_{c_1}$ there is a region of parameter space with three solutions, two of which are stable region and one unstable. At the critical value $\eps = \eps_{c_1}$, we identify a cubic fold point with location $\lambda\approx 0.9583$ (cf.~Fig.~\ref{fig:circle_side_eps_0_1_lam_0_6}). This example demonstrates the variability of solution multiplicity in the $(\eps,\lambda)$ parameter space.

\underline{Case $\eps = 0.1$:} At this fixed value, intervals of $\lambda$ in which multiple solutions of \eqref{eq:steady_problem} exist widens significantly. In general, we observe that as $\eps \rightarrow 0$, the potential multiplicity of solutions grows correspondingly. Indeed, the $\eps=0$ case discussed in Sec.~\ref{sec:conv_study} shows an infinite number of fold point centered on the values $\lambda = 4/9$. Figure \ref{fig:circle_side_eps_0_1_lam_0_6} (bottom row) shows three admitted solutions for $\eps = 0.1$ and $\lambda = 0.6$ with the position on the bifurcation curve highlighted. As $\eps$ decreases, we determine a second critical value $\eps_{c_2}\approx 0.0225$ where potential solution multiplicity increases to five.

\underline{Case $\eps = 0.01$:} In this parameter regime, the sharp features in equation \eqref{eq:steady_problem} become much more challenging to resolve and adaptive path-tracking is essential. The rapid spatial change in the solution occurs initially at the center of the domain where $\|\bu\|_{\infty} \approx 1-\eps$ followed by propagation radially outwards. To more accurately represent the solution in these regions, the adaptation approach outlined in Section \ref{sec:adapt_href} increases refinement in regions of large spatial gradient (cf.~Fig~\ref{fig:circle_tilted_eps_0_01_lam_0_444}). In panel 4 of Figure \ref{fig:circle_tilted_eps_0_01_lam_0_444}, the solution is rapidly changing in the center of the disk and the mesh is significantly refined in that region. 

\begin{figure}[t!]
\centering
  \resizebox{0.85\textwidth}{!}{\input{_py/square_all_bif_curves.tikz}}
  \resizebox{0.85\textwidth}{!}{\begin{tikzpicture}
\begin{groupplot} [
group style={group size = 3 by 2, horizontal sep = 0.05cm, vertical sep = 0.2cm},
title style={at={(current bounding box.north)}, anchor=south}]
\nextgroupplot[axis equal image, width=0.45\textwidth, xtick={-0, 1}, ytick={}, title={1}, xticklabels={}, yticklabels={}, ylabel={$\varepsilon = 0.3$}, xlabel={}, xmin=0, xmax=1, ymin=-1, ymax=0]
\addplot []
graphics [xmin=0,xmax=1,ymin=-1,ymax=0] { 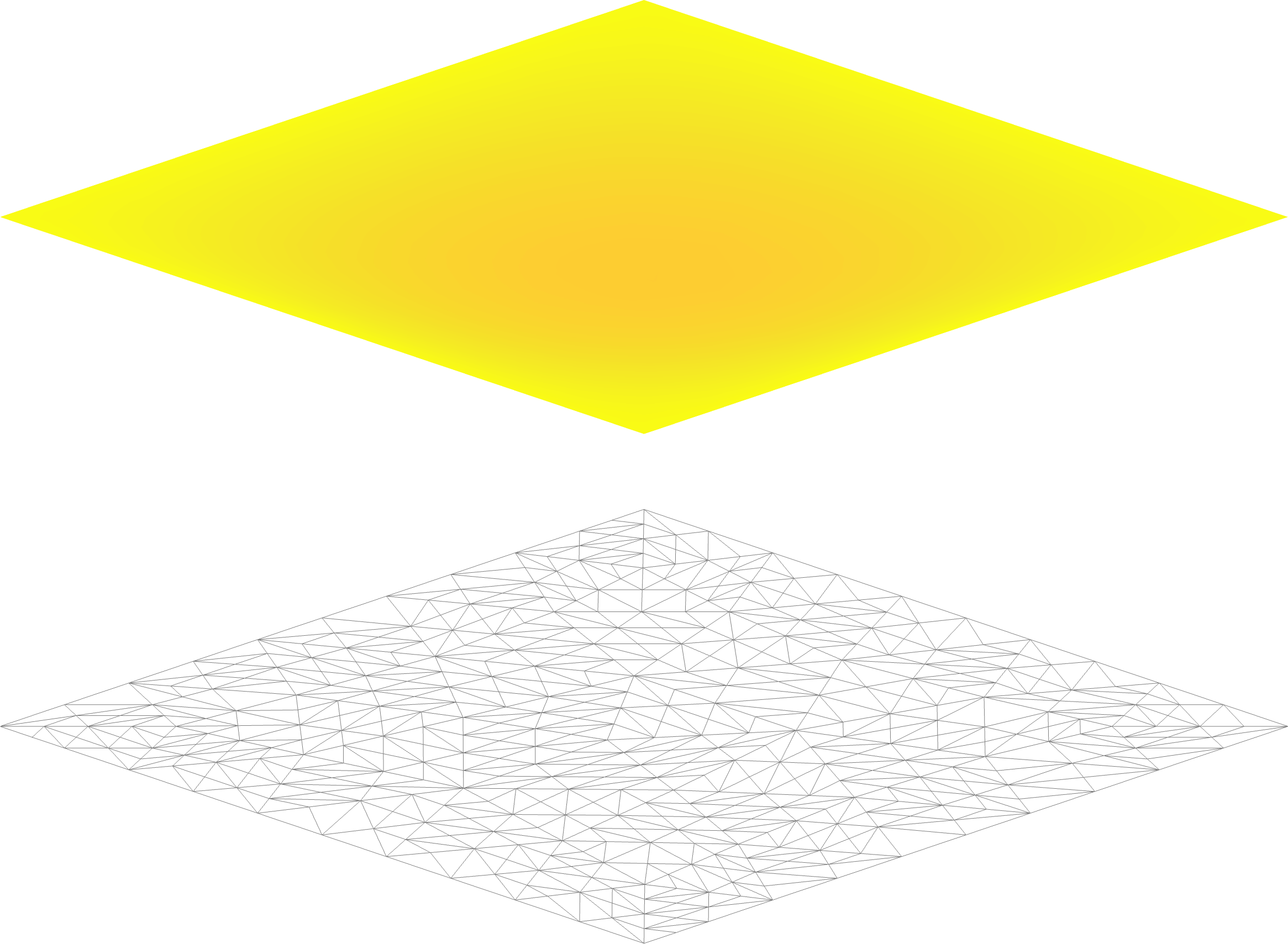};

\nextgroupplot[axis equal image, width=0.45\textwidth, xtick={}, ytick={7}, title={2}, xlabel={}, yticklabels={}, xticklabels={}, xmin=0, xmax=1, ymin=-1, ymax=0]
\addplot []
graphics [xmin=0,xmax=1,ymin=-1,ymax=0] { 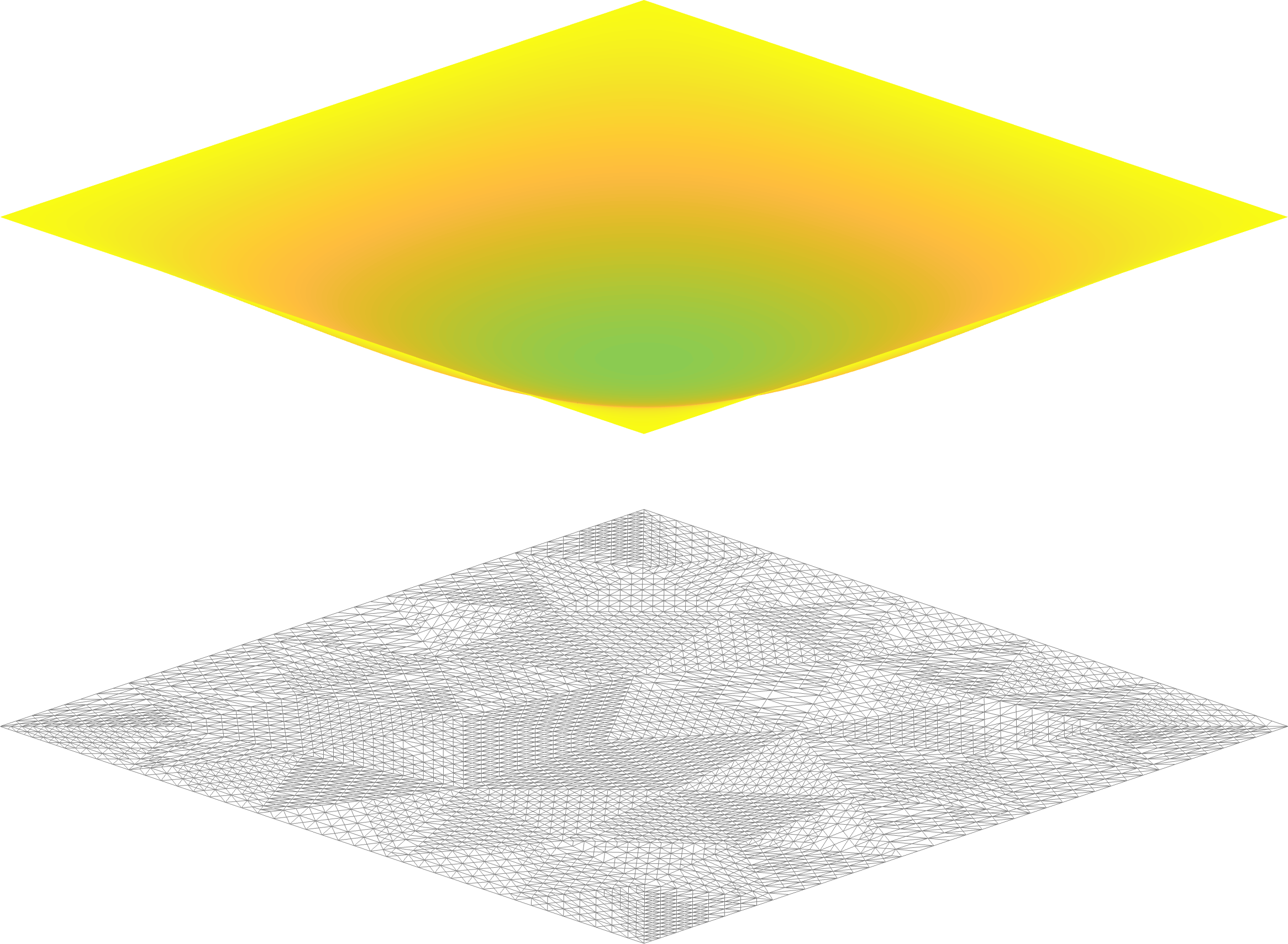};

\nextgroupplot[axis equal image, width=0.45\textwidth, xtick={}, ytick={}, title={3}, xlabel={}, xticklabels={}, yticklabels={}, xmin=0, xmax=1, ymin=-1, ymax=0]
\addplot []
graphics [xmin=0,xmax=1,ymin=-1,ymax=0] { 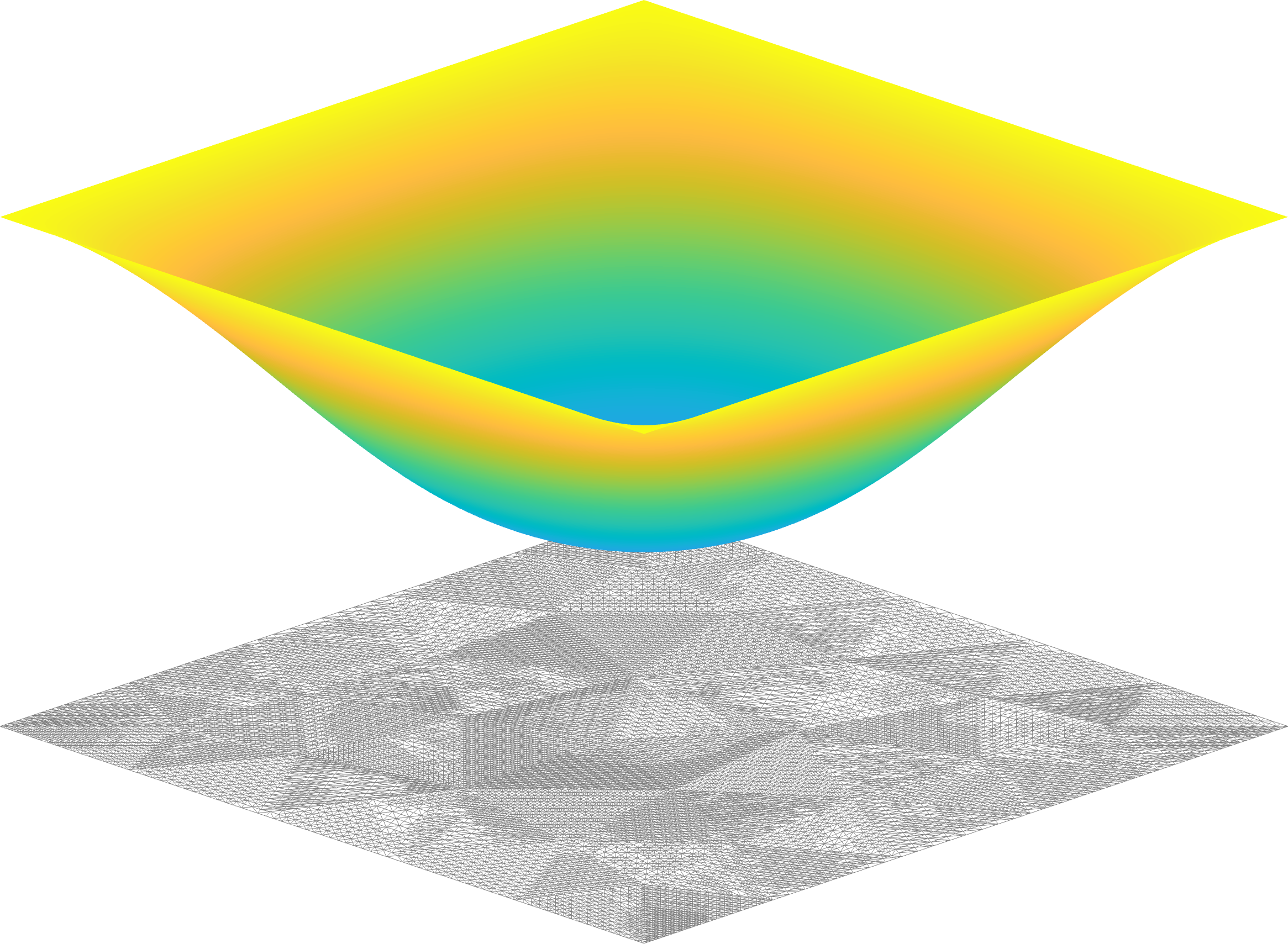};

\nextgroupplot[axis equal image, width=0.45\textwidth, xtick={}, ytick={}, xticklabels={}, yticklabels={}, ylabel={$\varepsilon = 0.1$}, xlabel={}, xmin=0, xmax=1, ymin=-1, ymax=0]
\addplot []
graphics [xmin=0,xmax=1,ymin=-1,ymax=0] { 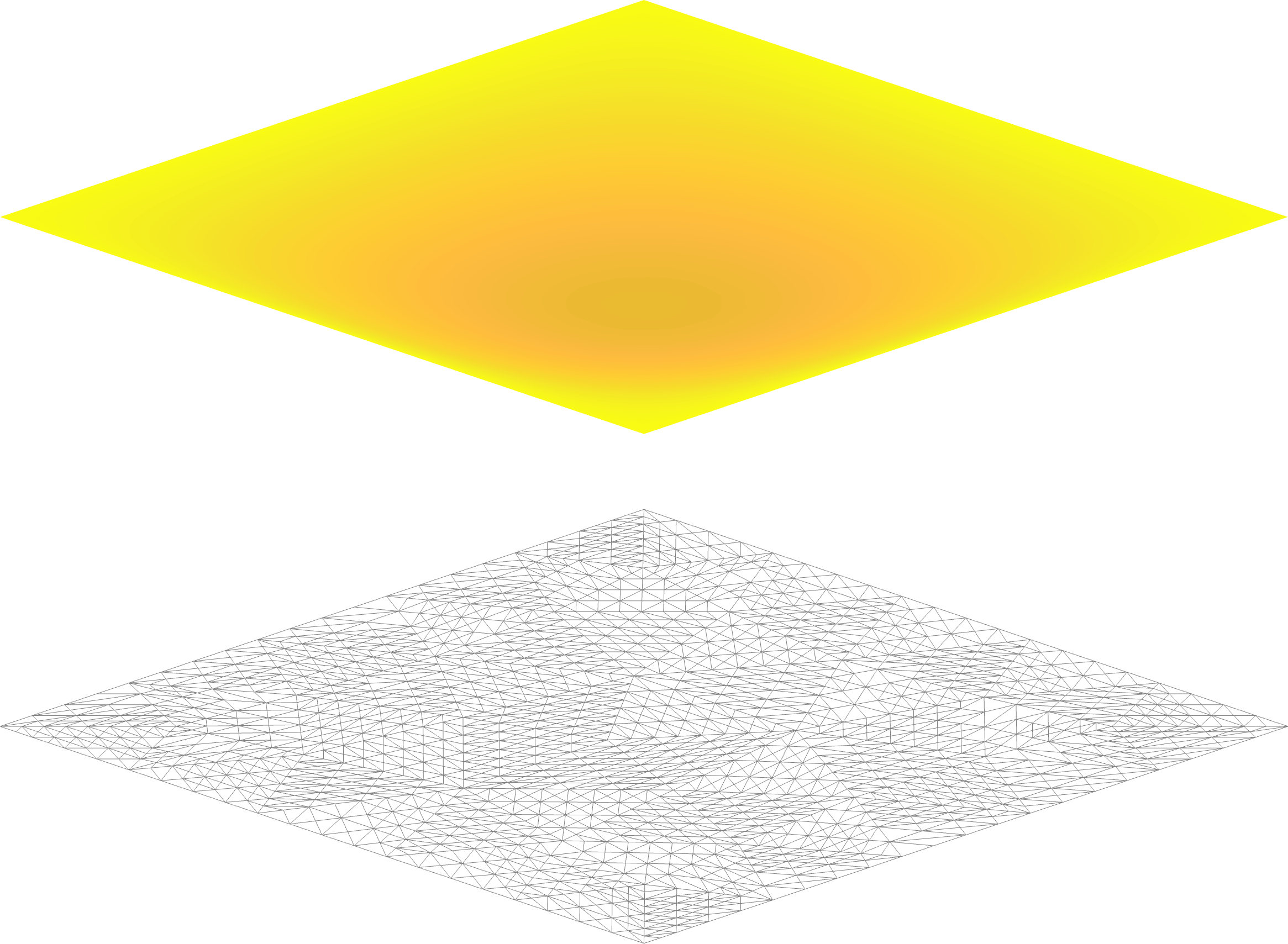};

\nextgroupplot[axis equal image, width=0.45\textwidth, xtick={}, ytick={}, xlabel={}, yticklabels={}, xticklabels={}, xmin=0, xmax=1, ymin=-1, ymax=0]
\addplot []
graphics [xmin=0,xmax=1,ymin=-1,ymax=0] { 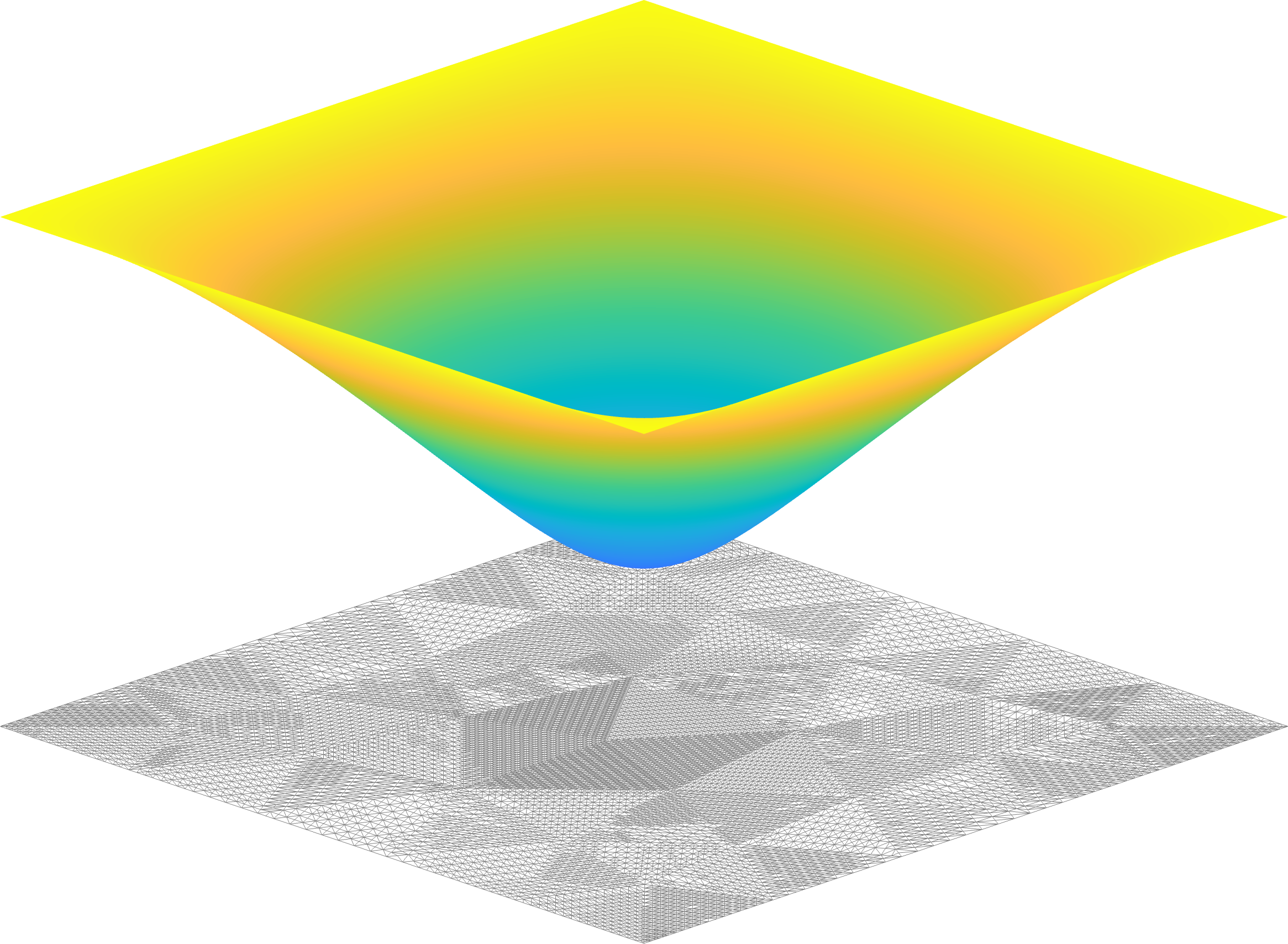};

\nextgroupplot[axis equal image, width=0.45\textwidth, xtick={}, ytick={}, xlabel={}, xticklabels={}, yticklabels={}, xmin=0, xmax=1, ymin=-1, ymax=0]
\addplot []
graphics [xmin=0,xmax=1,ymin=-1,ymax=0] { 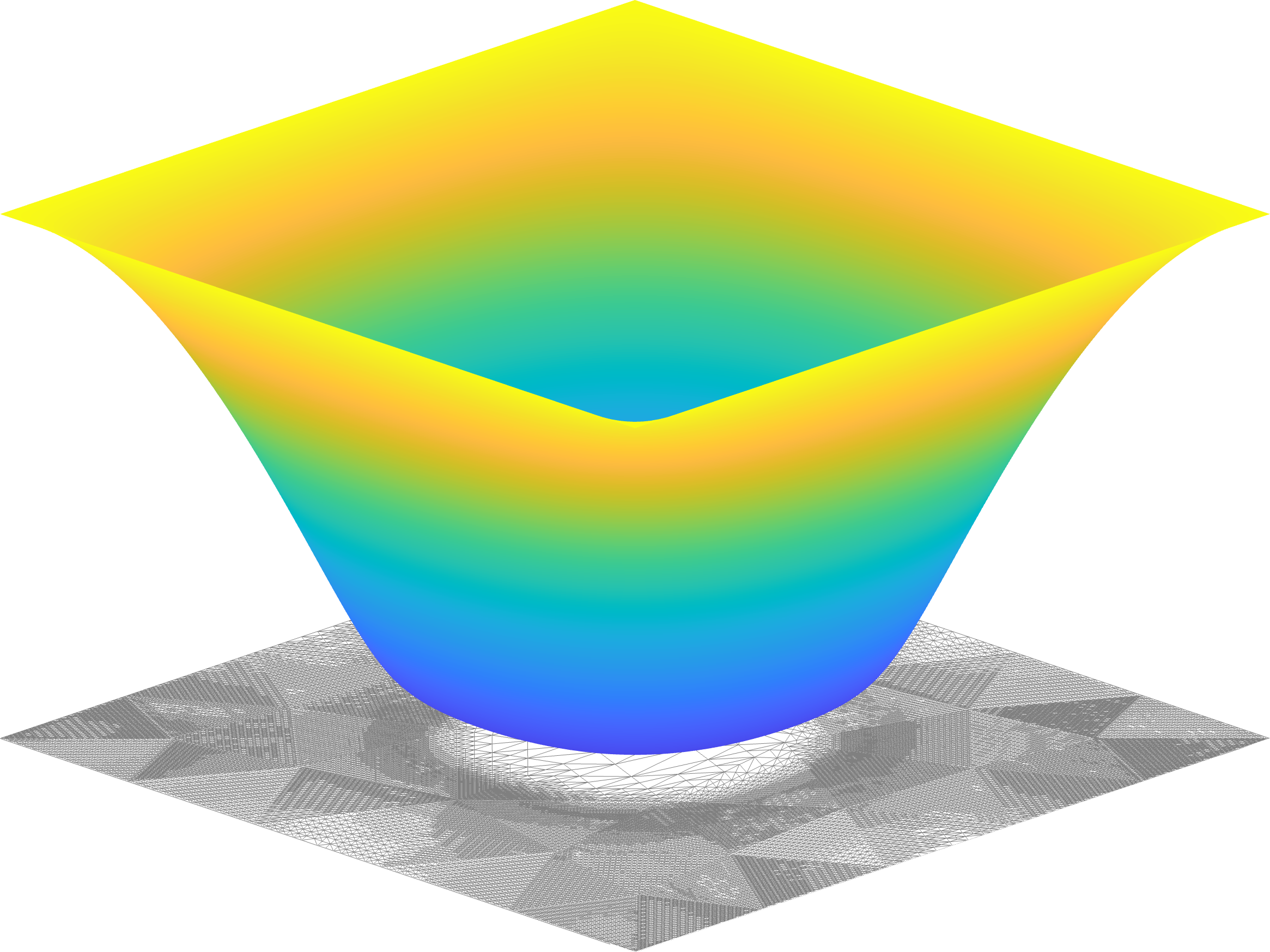};

\end{groupplot}\end{tikzpicture}}
 \colorbarMatlabParula{-1}{-0.75}{-0.5}{-0.25}{0}
 \caption{Solutions to \eqref{eq:steady_problem} on the unit square with $m=4$. Top row: Bifurcation curves of $\lambda$ against $\|\bu\|_2$ (\ref{line:bif_curve_eps_0.0225}) and $\|\bu\|_{\infty}$ (\ref{line:bif_curve_eps_0.1}) for a range of $\eps$ values \revtwo{where each curve represents a distinct family of solutions for a given value of $\eps$.} Middle row: Solution profiles $u$ and adapted mesh for points 1,2,3 (\ref{line:admitted_sol_eps_0_3}) on the $\eps = 0.3$ bifurcation curve. Bottom row: Solution profiles $u$ and adapted mesh for the points 1,2,3 (\ref{line:admitted_sol_eps_0_3}) on the $\eps = 0.1$ curve and $\lambda = 2$. }
 \label{fig:square_bif_curves_eps_0_3_to_0_1}
\end{figure}

In Fig.~\ref{fig:circle_tilted_eps_0_01_lam_0_444} (panel 5), the region of large gradient in the solution has propagated outwards and the mesh is refined in the vicinity of this feature. In the center of the domain, where the solution is largely constant, the mesh has been significantly coarsened. The locations of the five solution plots for the fixed value $\lambda = 0.444$ are shown in Fig.~\ref{fig:circle_tilted_eps_0_01_lam_0_444} along the bifurcation curve (\ref{line:bif_curve_eps_0.01_lam_0.444}) together with an inset showing how the multiplicity changes in a very small $\lambda$ interval.

\begin{figure}[t!]
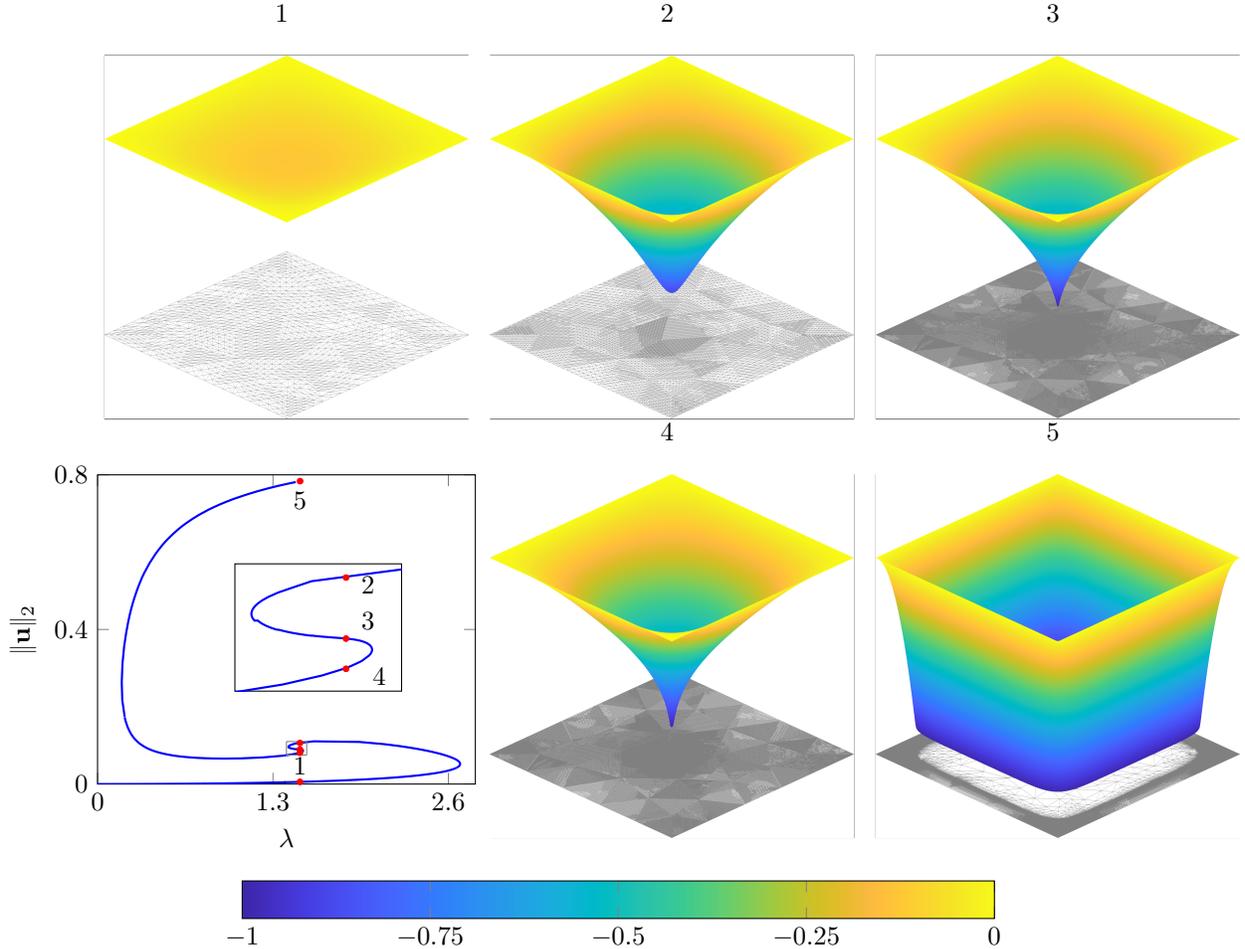

\centering
 \raisebox{-0.5\height}{\input{_py/square_tilted_eps_0_01_lam_1_5.tikz}}
 \colorbarMatlabParula{-1}{-0.75}{-0.5}{-0.25}{0}
\begin{picture}(0,0)
\put(-300, 100){\input{_py/square_tilted_eps_0_01_lam_0_444_inset.tikz}}
\end{picture}
 \caption{Solution of \eqref{eq:steady_problem} on the unit square for $\eps = 0.01$ with five profiles $u$ and adapted meshes for $\lambda \approx 1.5$. The plotted solutions (panels 1-5) are shown in their positions (\ref{line:point_eps_0.01_lam_0.444_1}) along the bifurcation curve (\ref{line:bif_curve_eps_0.01_lam_0.444}).}
 \label{fig:square_tilted_eps_0_01_lam}
\end{figure}

We observe that as $\eps \rightarrow 0$, the multiplicity of solutions begins to approach infinity at a fixed value of $\lambda = \lambda_c$. The radially-symmetric solution developed in Section \ref{sec:conv_study} shows that the bifurcation appears as a spiral where an arbitrary number of solutions exist in the neighborhood of $\lambda_c = 4/9$. In our numerical results, we observe a close approach to this spiral point, e.g. in the case $\eps = 0.01$ we observe five solutions for certain intervals of $\lambda$. As the multiplicity of solutions increase, the tightness of the fold points on the bifurcation curve significantly increases. This makes traversing the bifurcation curve particularly challenging in the vicinity of $\lambda_c$ and necessitates choosing a very small arc-length step-size $ds$ in such regions of the curve.

\subsection{Regularized problem on square domain} \label{sec:reg_unit_square}

In this section we present results for the regularized problem \eqref{eq:steady_problem} on the square domain, $\Omega = \{ (x_1,x_2) \in \mathbb{R}^2 \, \mid \, 0 \leq x_d \leq 1, d = 1,2\}$ and observe similar solution multiplicity structures as $\eps$ decreases (cf.~Fig.~\ref{fig:square_bif_curves_eps_0_3_to_0_1}).

\underline{Case $\eps = 0.3$:} In Fig.~\ref{fig:square_bif_curves_eps_0_3_to_0_1} (middle row) we show solutions of \eqref{eq:steady_problem} for $\eps = 0.3$, $m=4$ and $\lambda = [1.5, 3, 4.5]$. We observe a unique solution for each $\lambda>0$. Similar to observed solutions in the circular case, the deflection is initially largest in the center until the source terms balance and there is a flat central region with increasing deflection outwards. The adaptive refinement procedure is able to add resolution to areas of the domain with sharp gradient and coarsen where the solution has low gradient (constant solution).

\underline{Case $\eps = 0.1$}:
As $\eps$ decreases, we observe increased solution multiplicity in \eqref{eq:steady_problem} and solutions with large spatial gradients. We identify a critical value $\eps = \eps_{c_1} \approx 0.25966$ at which the multiplicity of solutions increases through a cubic fold point (cf.~Fig.~\ref{fig:square_bif_curves_eps_0_3_to_0_1}). For $\eps < \eps_{c_1}$, we observe intervals in $\lambda$ with three solutions (two stable regions, one unstable). At the value $\eps = 0.1$, we plot in Fig.~\ref{fig:square_bif_curves_eps_0_3_to_0_1} (bottom row) three solutions of \eqref{eq:steady_problem} for $\lambda = 2$

\underline{Case $\eps = 0.01$:} Decreasing $\eps$ further, we identify intervals of $\lambda$ where \eqref{eq:steady_problem} admits five solutions and observe the formation of a spiral in the bifurcation diagram (cf.~Fig.~\ref{fig:square_tilted_eps_0_01_lam}). The numerical method exhibits significant robustness in its ability to approximate the solution for the whole bifurcation curve, particularly in tight regions where solution multiplicity changes over very small intervals of $\lambda$. In Figure \ref{fig:square_tilted_eps_0_01_lam} we display five admitted solutions for $\lambda = 1.5$ and note that mesh adaptation is essential for their resolution.

\subsection{Two-dimensional annulus for symmetry breaking} \label{sec:symm_break}
In this section, we present computed solutions of \eqref{eq:steady_problem} for the annulus domain $\Omega = \{ \bx \in \mathbb{R}^2 \, \mid \, r_1 \leq |\bx| \leq 1 \}$. This example explores changes in solution structure due to perturbations in domain topology, which MEMS can be engineered with (cf.~Fig.\ref{fig:intro_diagram}). In particular, we compare with and expand upon results from \cite{Pelesko_sym_breaking}. Specifically, we use tools presented in Section \ref{sec:branch_switch} to traverse the radially symmetric solution branch until symmetry-breaking bifurcations are \revSec{detected by zeroes in the determinant} of the Jacobian, det$f_{\bu}^T (\bu^{(k)}, \lambda^{(k)})$. At such points, we solve the system of equations \eqref{eq:branch_point_locator} to locate the branch point and then travel along the new asymmetric solution branch. 

For simplicity of exhibition, we fix $r_1 =0.1$ and perform continuation in $\lambda$ while solving for the deflection $u$ over the domain. We apply localized mesh refinement for improved resolution in regions of the domain where solution gradients are large (as determined by the hierarchical error indicator). In our formulation the boundaries of the annulus are held fixed at $u = 0$ on $\partial \Omega$. Accordingly, deflection occurs in the middle of the ring first and increases radially outward and inward from there. 

\begin{figure}[t!]
\centering
 \raisebox{-0.0\height}{
    \input{_py/annulus_branch_bif_curve_eps_0.tikz}
 }
  \resizebox{0.43\textwidth}{!}{
    \begin{tikzpicture}
\begin{groupplot} [
group style={group size = 3 by 3, horizontal sep = 0.2cm, vertical sep = 0.75cm},
title style={at={(current bounding box.north)}, anchor=south}]
\nextgroupplot[axis equal image, axis lines=none, width=0.45\textwidth, xtick={0, 1}, ytick={0}, title={1}, xticklabels={}, yticklabels={}, ylabel={$u$}, xmin=0, xmax=1, ymin=-1, ymax=0]
\addplot []
graphics [xmin=0,xmax=1,ymin=-1,ymax=0] { 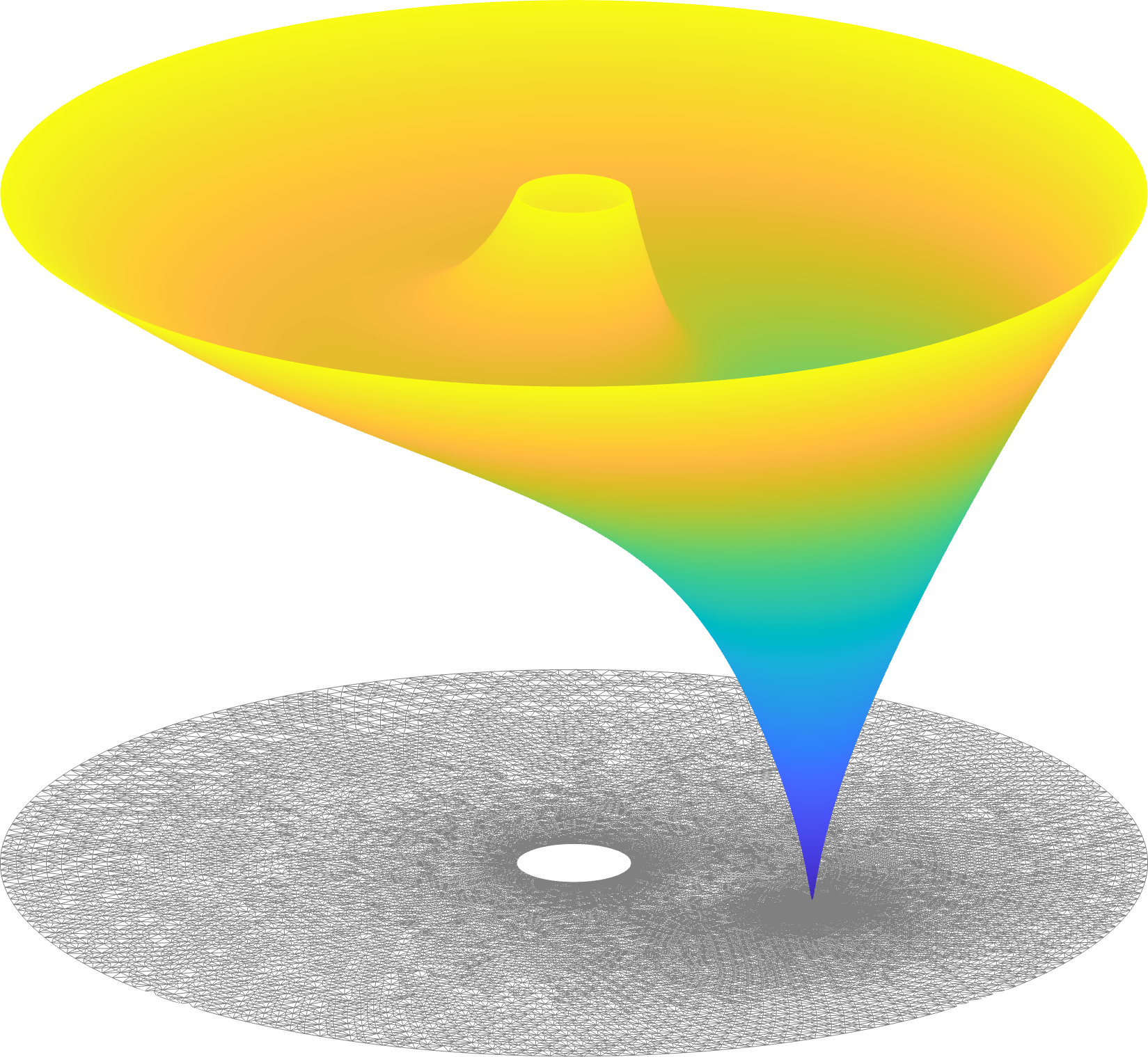};

\nextgroupplot[axis equal image, axis lines=none, width=0.45\textwidth, xtick={0, 1}, ytick={0}, title={2}, xticklabels={}, yticklabels={}, xmin=0, xmax=1, ymin=-1, ymax=0]
\addplot []
graphics [xmin=0,xmax=1,ymin=-1,ymax=0] { 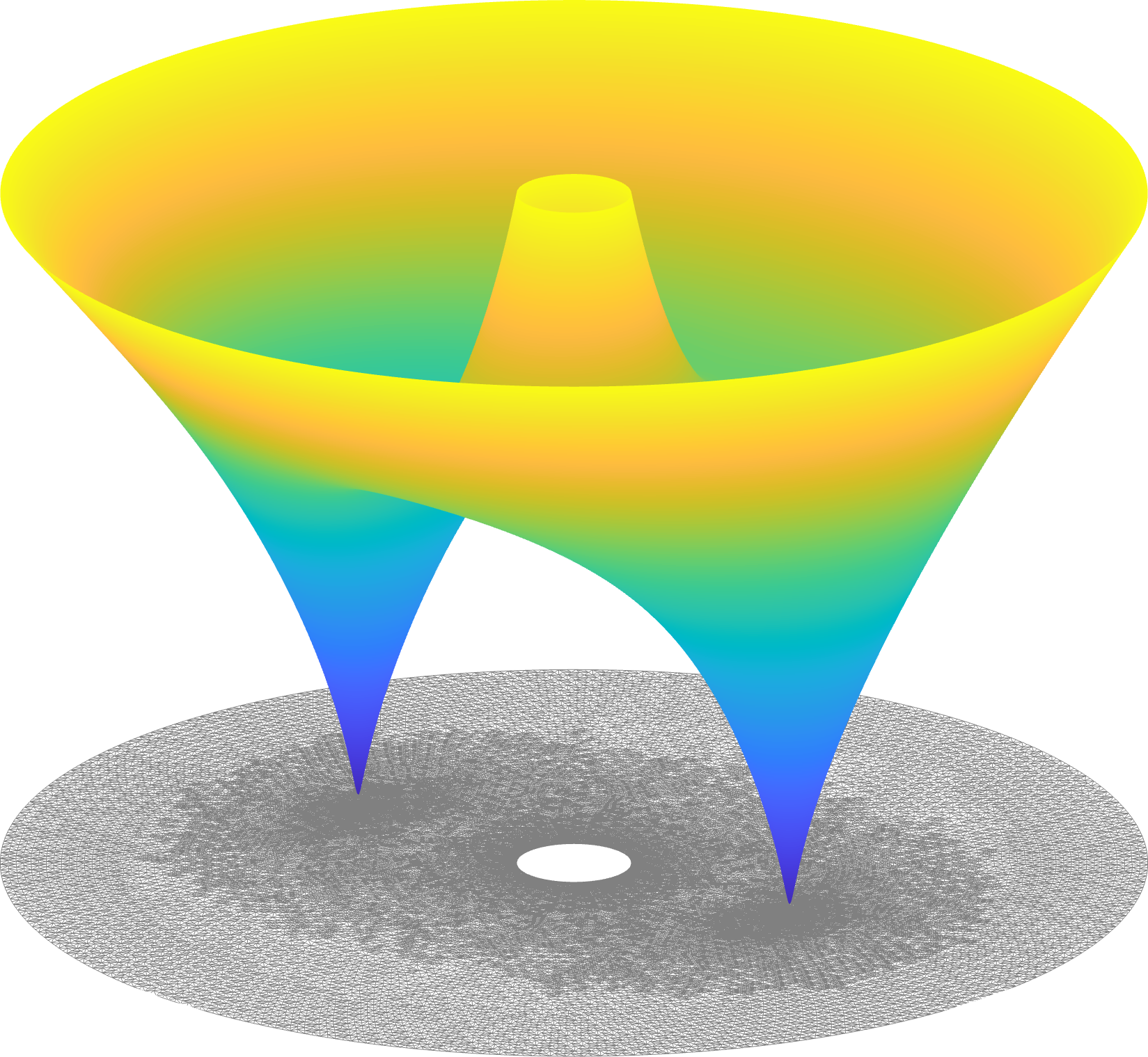};

\nextgroupplot[axis equal image, axis lines=none, width=0.45\textwidth, xtick={0, 1}, ytick={0}, title={3}, xticklabels={}, yticklabels={}, xmin=0, xmax=1, ymin=-1, ymax=0]
\addplot []
graphics [xmin=0,xmax=1,ymin=-1,ymax=0] { 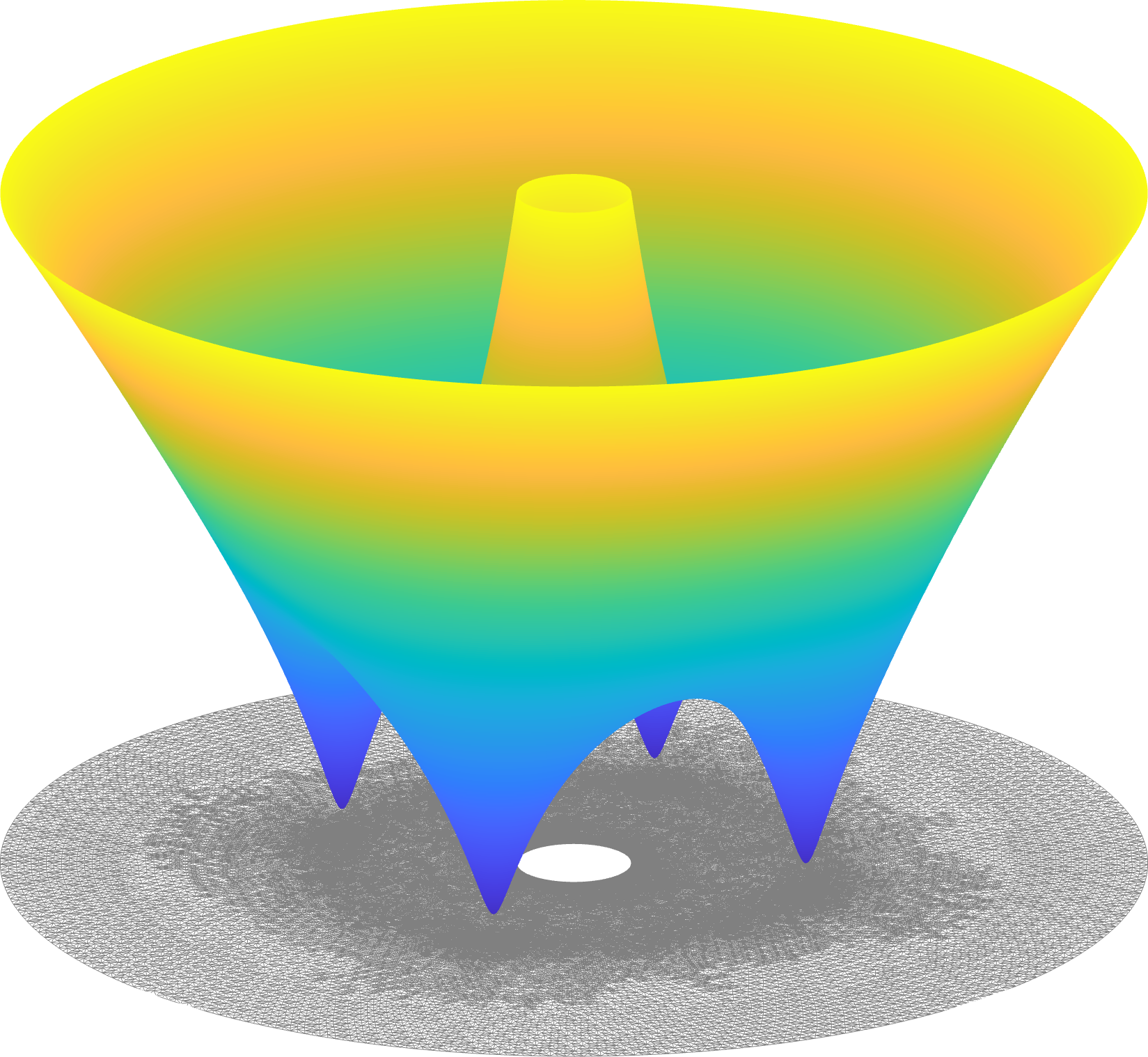};

\nextgroupplot[axis equal image, axis lines=none, width=0.45\textwidth, xtick={0, 1}, ytick={0}, title={4}, xticklabels={}, yticklabels={}, xlabel={$x$}, ylabel={$u$}, xmin=0, xmax=1, ymin=-1, ymax=0]
\addplot []
graphics [xmin=0,xmax=1,ymin=-1,ymax=0] { 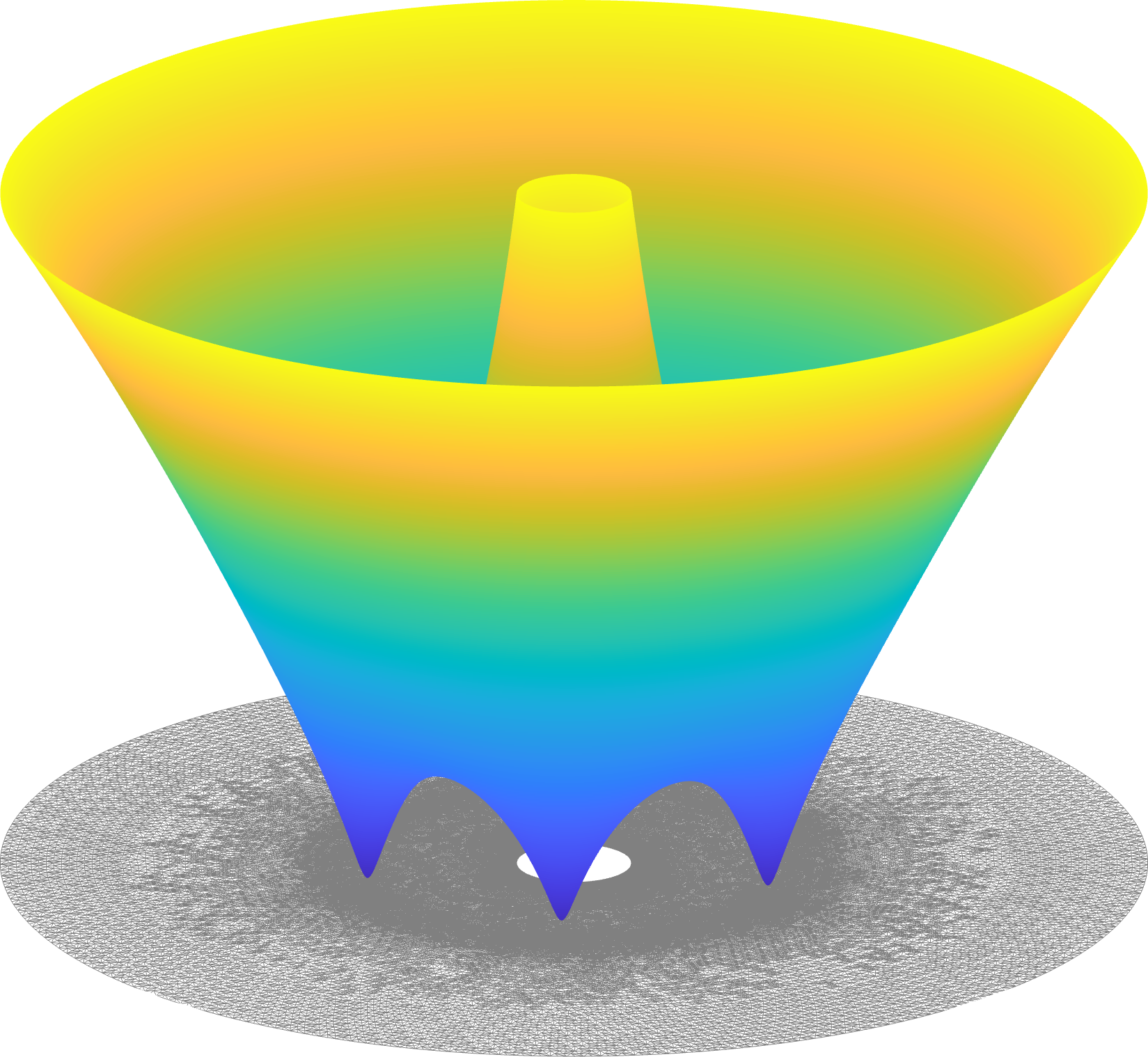};

\nextgroupplot[axis equal image, axis lines=none, width=0.45\textwidth, xtick={0, 1}, ytick={0}, title={5}, xticklabels={}, yticklabels={}, xlabel={$x$}, ylabel={$u$}, xmin=0, xmax=1, ymin=-1, ymax=0]
\addplot []
graphics [xmin=0,xmax=1,ymin=-1,ymax=0] { 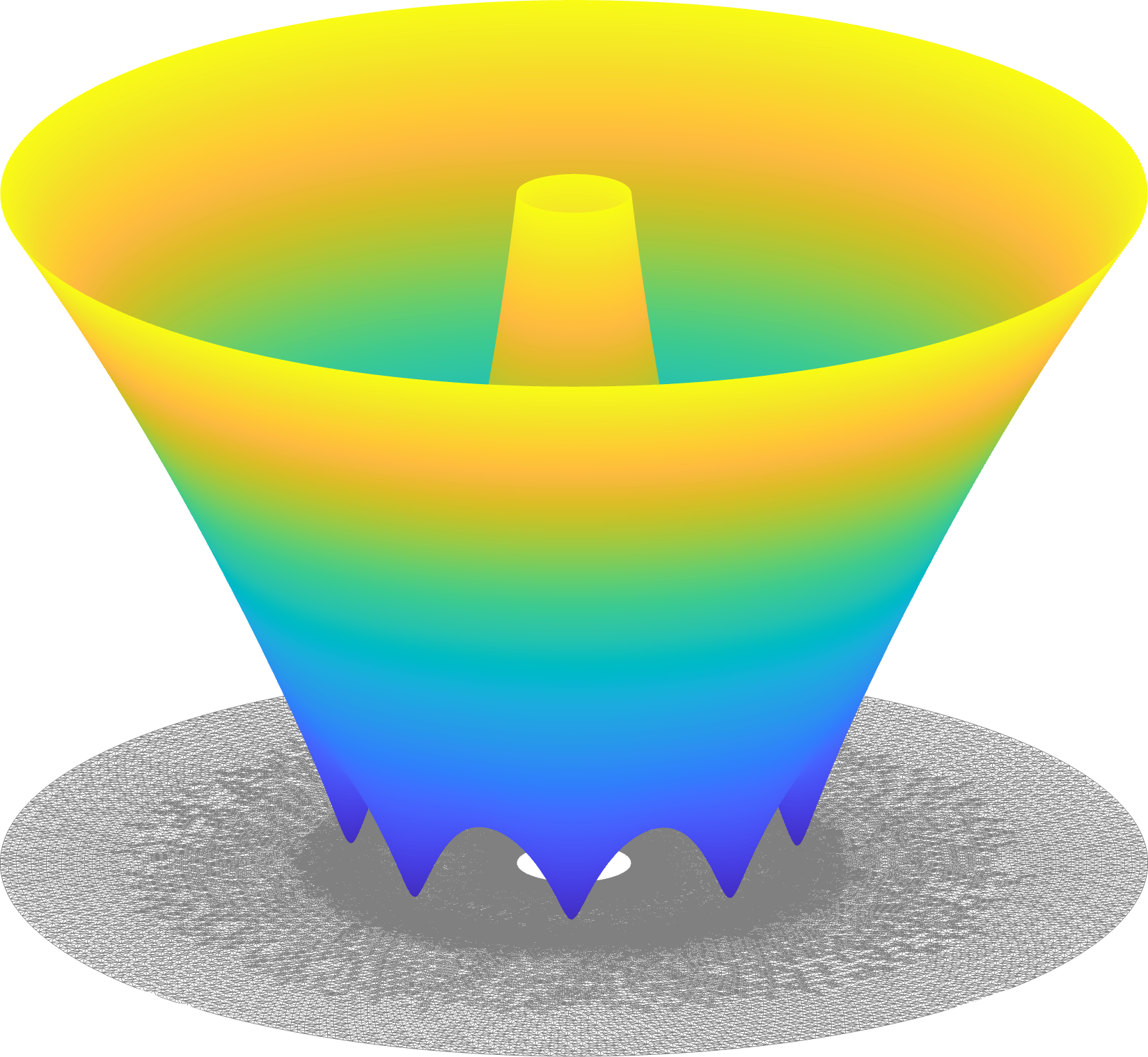};

\nextgroupplot[axis equal image, axis lines=none, width=0.45\textwidth, xtick={0, 1}, ytick={0}, title={6}, xticklabels={}, yticklabels={}, xlabel={$x$}, ylabel={$u$}, xmin=0, xmax=1, ymin=-1, ymax=0]
\addplot []
graphics [xmin=0,xmax=1,ymin=-1,ymax=0] { 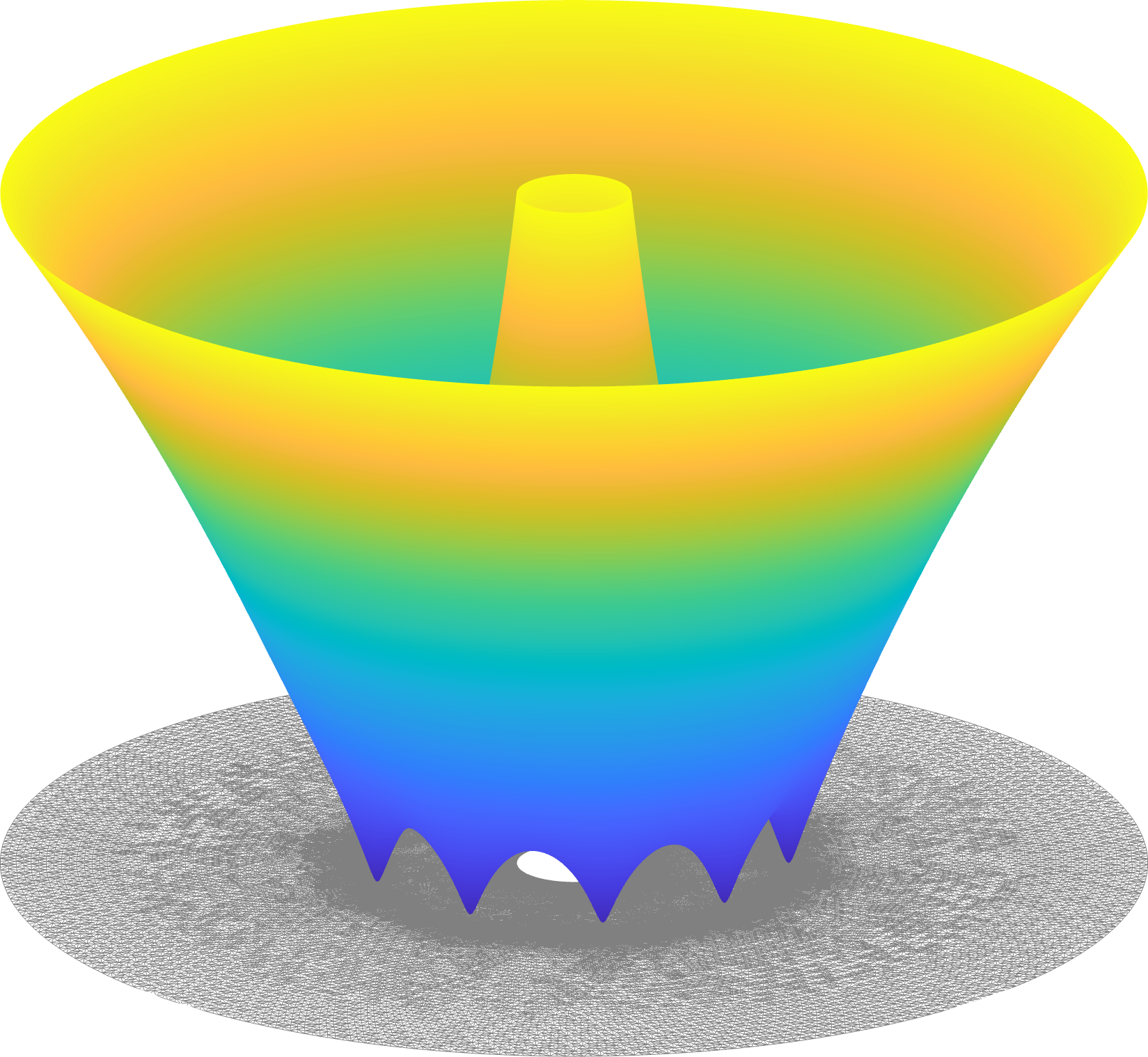};

\nextgroupplot[axis equal image, axis lines=none, width=0.45\textwidth, xtick={0, 1}, ytick={0}, title={7}, xticklabels={}, yticklabels={}, xlabel={$x$}, ylabel={$u$}, xmin=0, xmax=1, ymin=-1, ymax=0]
\addplot []
graphics [xmin=0,xmax=1,ymin=-1,ymax=0] { 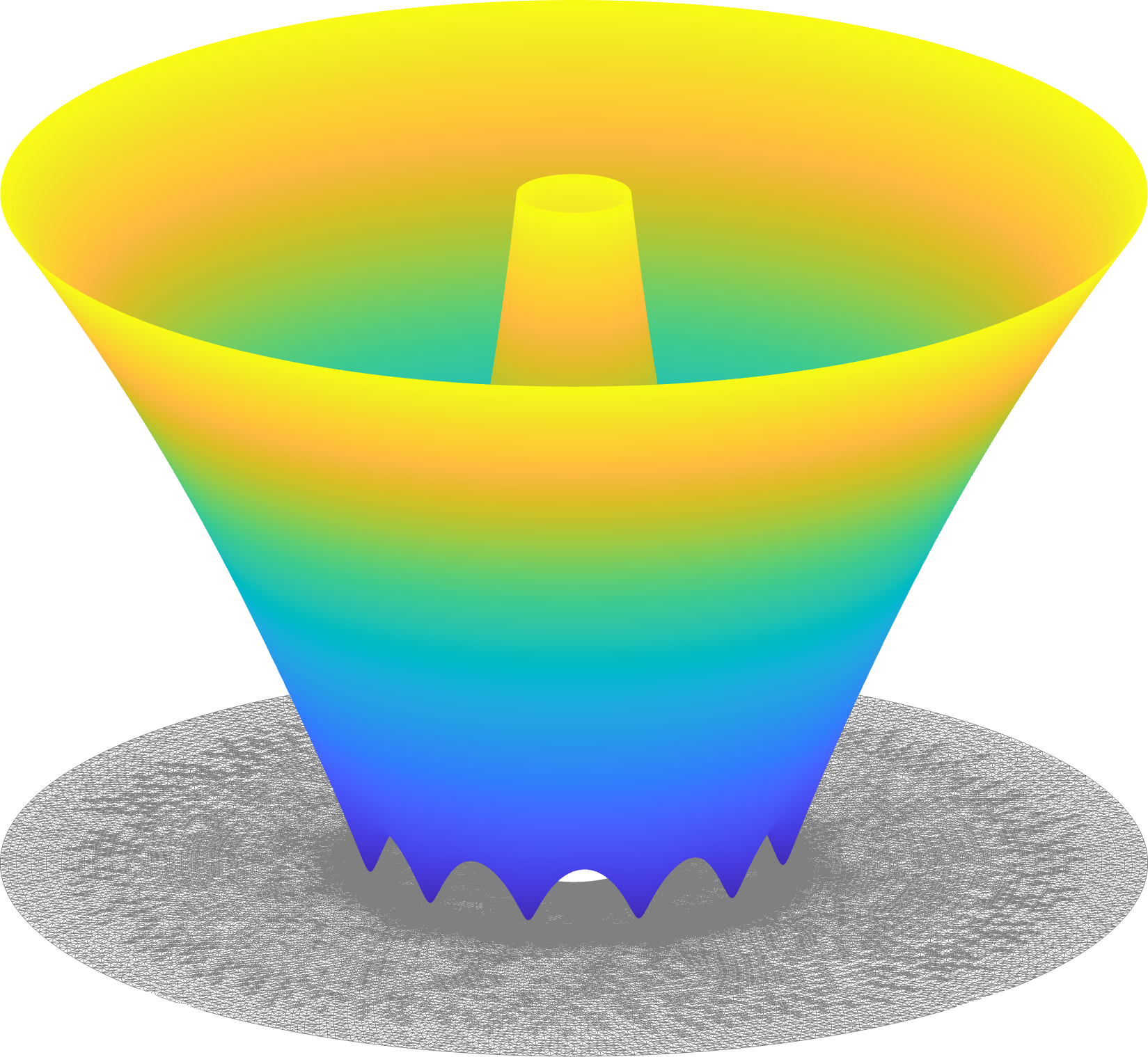};

\nextgroupplot[axis equal image, axis lines=none, width=0.45\textwidth, xtick={0, 1}, ytick={0}, title={8}, xticklabels={}, yticklabels={}, xlabel={$x$}, ylabel={$u$}, xmin=0, xmax=1, ymin=-1, ymax=0]
\addplot []
graphics [xmin=0,xmax=1,ymin=-1,ymax=0] { 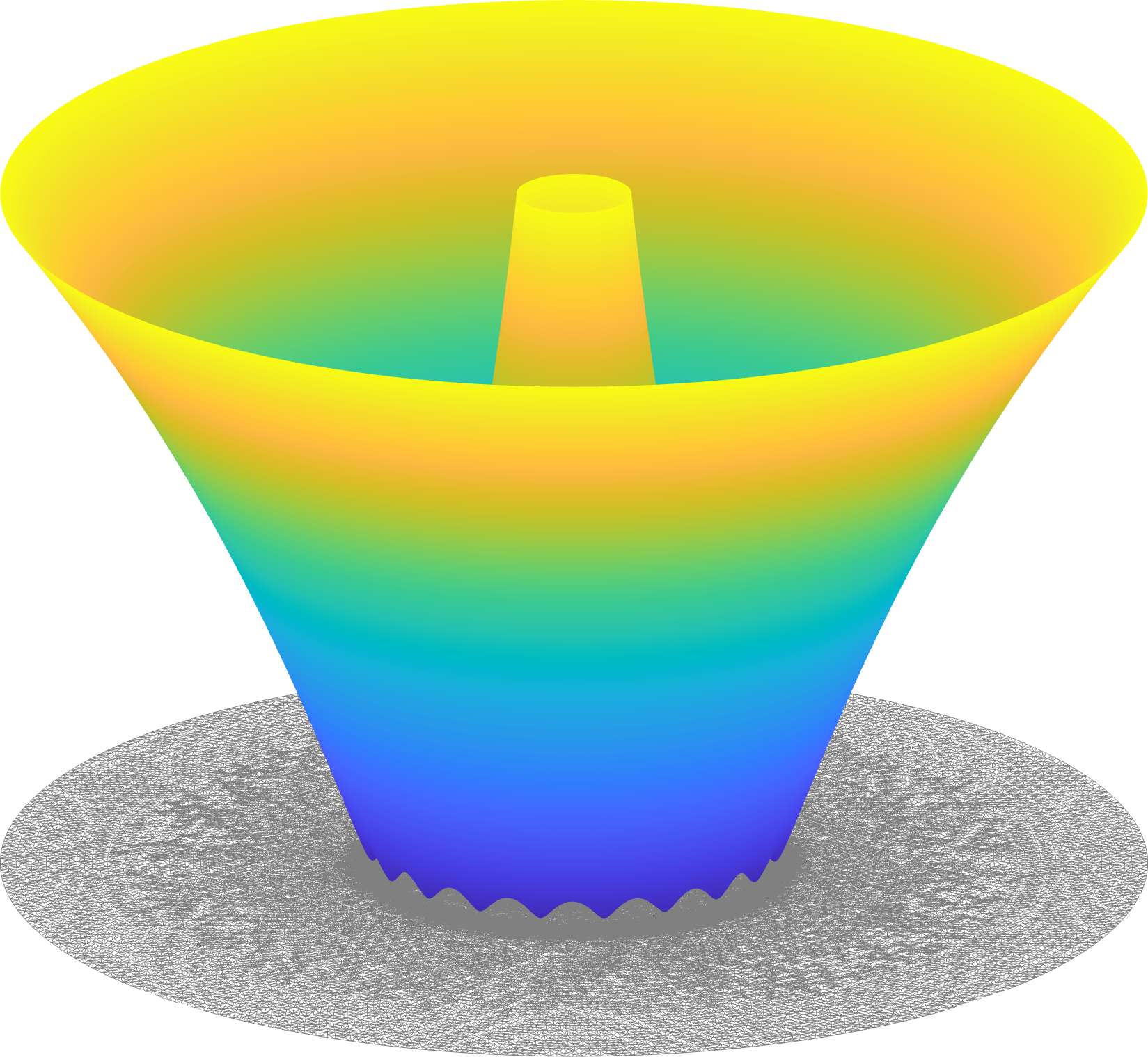};

\nextgroupplot[axis equal image, axis lines=none, width=0.45\textwidth, xtick={0, 1}, ytick={0}, title={9}, xticklabels={}, yticklabels={}, xlabel={$x$}, ylabel={$u$}, xmin=0, xmax=1, ymin=-1, ymax=0]
\addplot []
graphics [xmin=0,xmax=1,ymin=-1,ymax=0] { 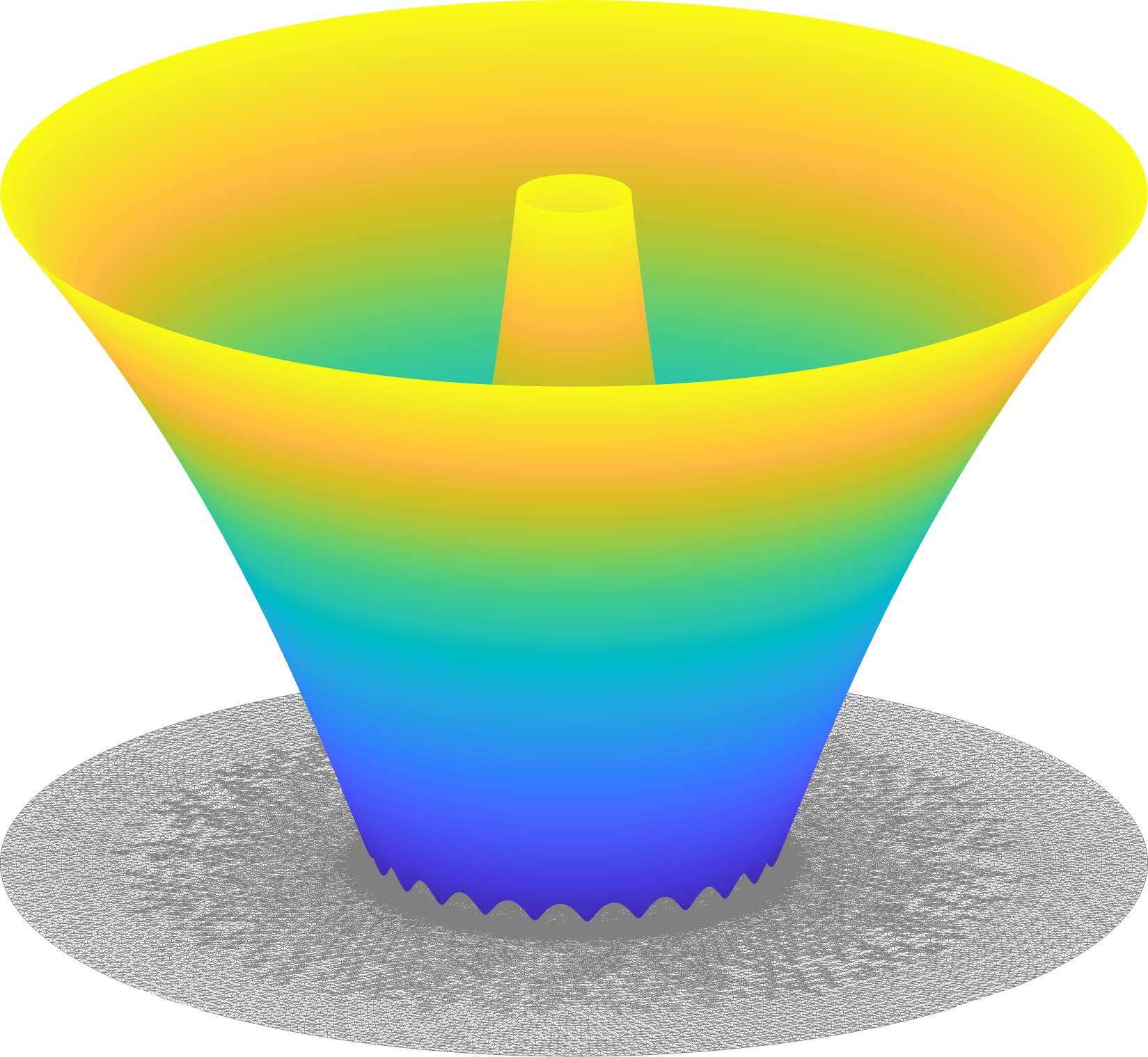};

\end{groupplot}\end{tikzpicture}
    }
\colorbarMatlabParula{-1}{-0.75}{-0.5}{-0.25}{0}
 \caption{Solution of \eqref{eq:steady_problem} for $\eps=0$ on the annular region $0.1<|\bx|<1$. Left: Bifurcation curves for $\|\bu\|_2$ (solid blue)(\ref{line:bif_curve_annulus_branch_1}) and $\|\bu\|_\infty$ (solid green) (\ref{line:bif_curve_annulus_branch_1_uinf}) showing main branch and first nine asymmetric branches. Right: Solutions and associated computational meshes for each of the labeled points.}
 \label{fig:annulus_bif_curves}
\end{figure}

In Fig.~\ref{fig:annulus_bif_curves} we show bifurcation curves for both $\|\bu\|_2$ (\ref{line:bif_curve_annulus_branch_1}) and $\|\bu\|_\infty$ (\ref{line:bif_curve_annulus_branch_1_uinf}) for $\eps=0$ computed using our adaptive path-tracking method. Bifurcation curves are characterized by a radially symmetric branch with asymmetric branches that are characterized by their minima. Previous studies \cite{Feng2005} have theoretically predicted the existence of an infinite number of branching points and here we numerically resolve the first nine of these branches. 
 In Fig.~\ref{fig:annulus_bif_curves}, we also plot the final computed solution along each asymmetric branch together with the computational mesh.

\begin{figure}
\centering
 \raisebox{-0.0\height}{
    \input{_py/annulus_branch_bif_curve_eps_gr_0.tikz}
 }
  \resizebox{0.43\textwidth}{!}{
    \begin{tikzpicture}
\begin{groupplot} [
group style={group size = 3 by 3, horizontal sep = 0.2cm, vertical sep = 0.75cm},
title style={at={(current bounding box.north)}, anchor=south}]
\nextgroupplot[axis equal image, axis lines=none, width=0.45\textwidth, xtick={0, 1}, ytick={0}, title={1}, xticklabels={}, yticklabels={}, ylabel={$\varepsilon = 0.2$}, xmin=0, xmax=1, ymin=-1, ymax=0]
\addplot []
graphics [xmin=0,xmax=1,ymin=-1,ymax=0] { 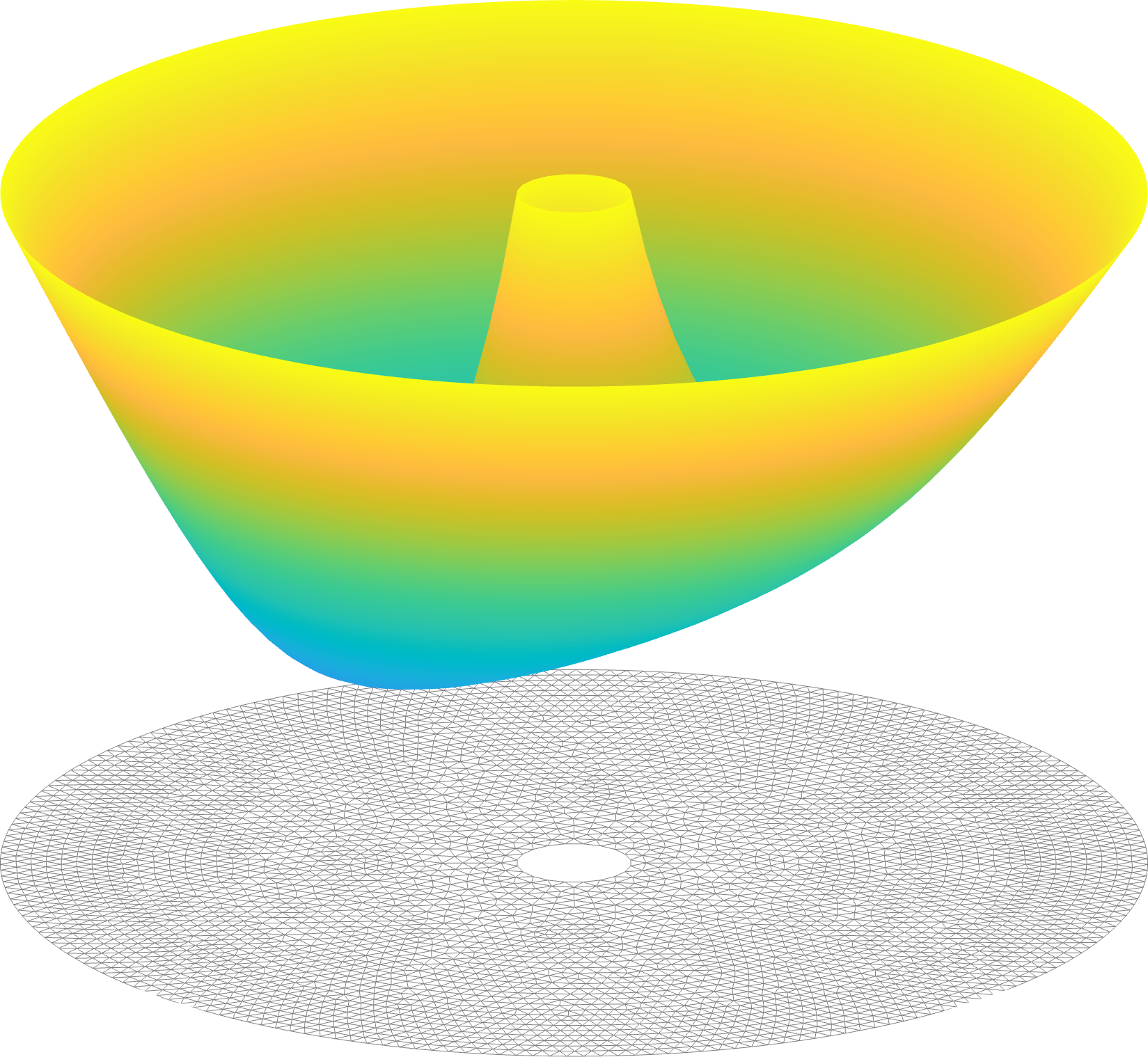};

\nextgroupplot[axis equal image, axis lines=none, width=0.45\textwidth, xtick={0, 1}, ytick={0}, title={2}, xticklabels={}, yticklabels={}, xmin=0, xmax=1, ymin=-1, ymax=0]
\addplot []
graphics [xmin=0,xmax=1,ymin=-1,ymax=0] { 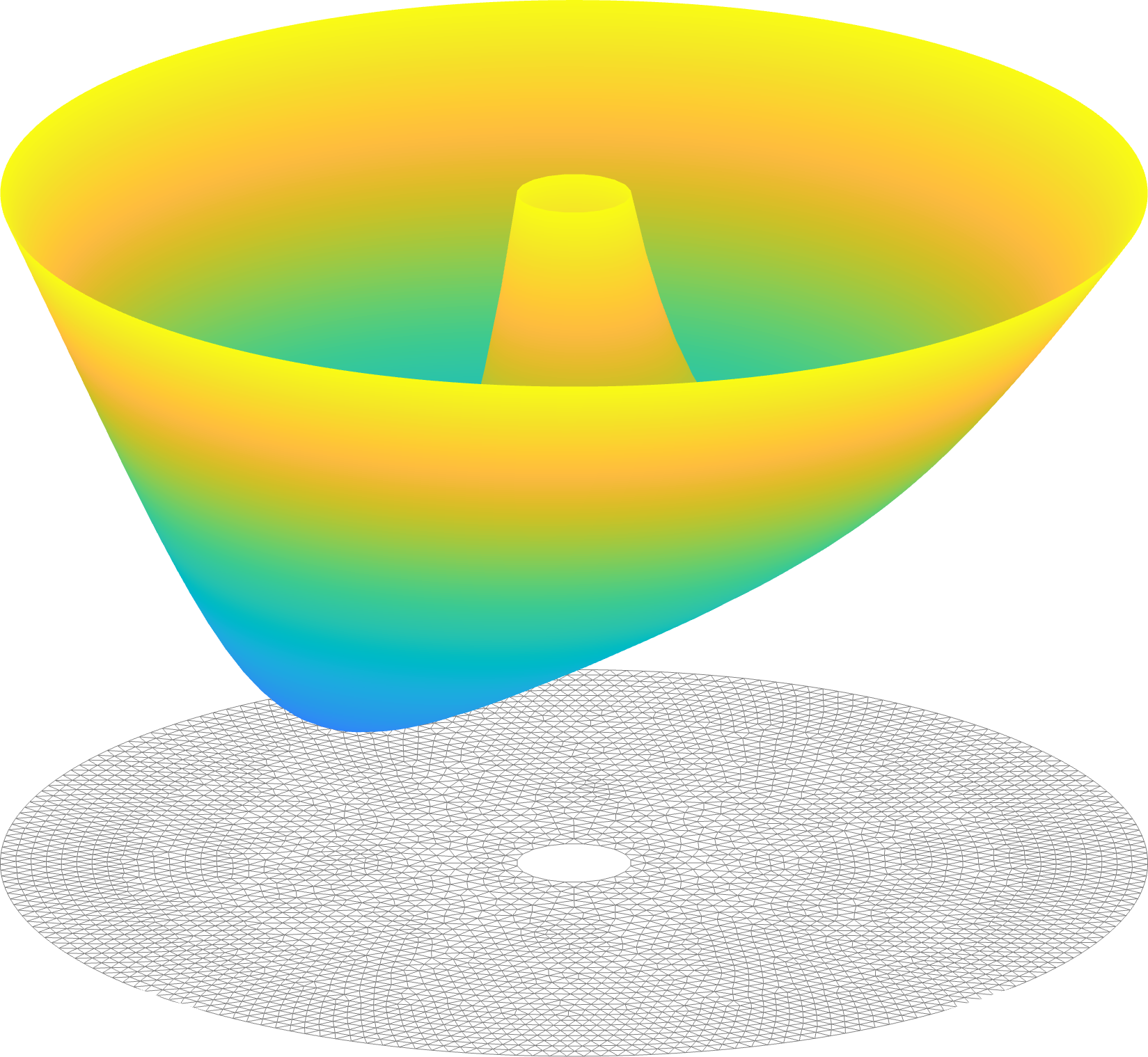};

\nextgroupplot[axis equal image, axis lines=none, width=0.45\textwidth, xtick={0, 1}, ytick={0}, title={3}, xticklabels={}, yticklabels={}, xmin=0, xmax=1, ymin=-1, ymax=0]
\addplot []
graphics [xmin=0,xmax=1,ymin=-1,ymax=0] { 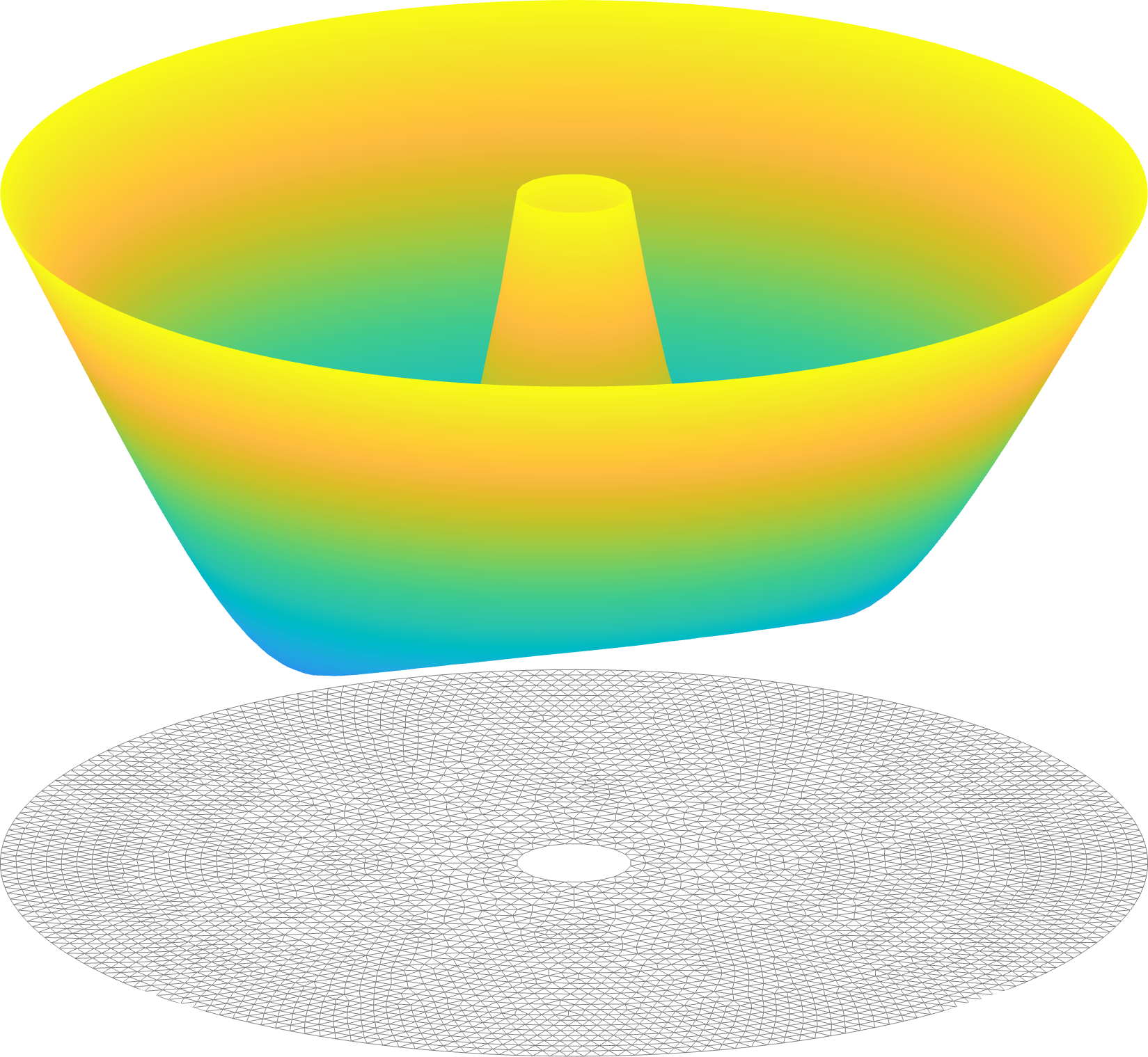};

\nextgroupplot[axis equal image, axis lines=none, width=0.45\textwidth, xtick={0, 1}, ytick={0}, title={4}, xticklabels={}, yticklabels={}, ylabel={$\varepsilon = 0.1$}, xmin=0, xmax=1, ymin=-1, ymax=0]
\addplot []
graphics [xmin=0,xmax=1,ymin=-1,ymax=0] { 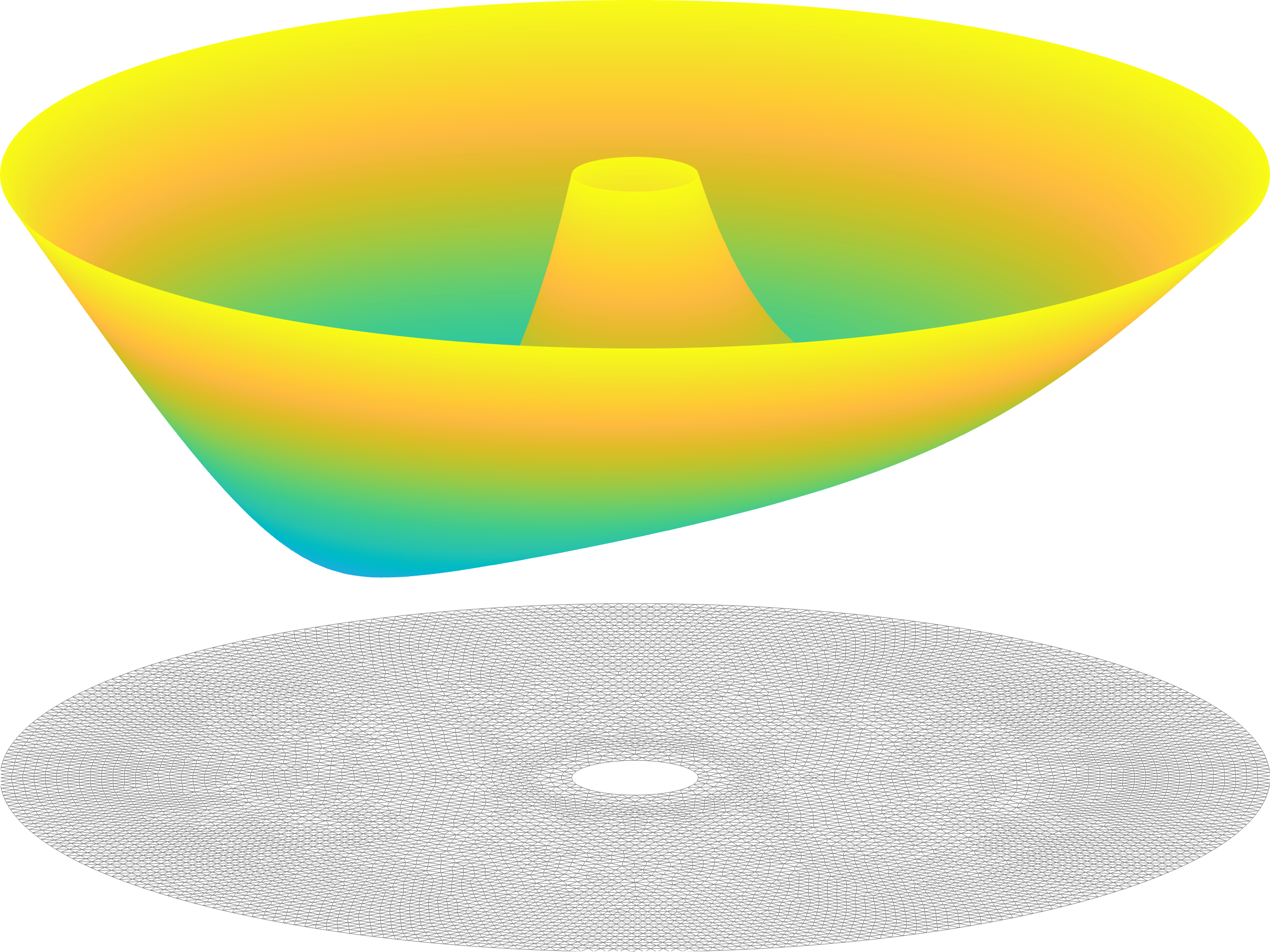};

\nextgroupplot[axis equal image, axis lines=none, width=0.45\textwidth, xtick={0, 1}, ytick={0}, title={5}, xticklabels={}, yticklabels={}, xmin=0, xmax=1, ymin=-1, ymax=0]
\addplot []
graphics [xmin=0,xmax=1,ymin=-1,ymax=0] { 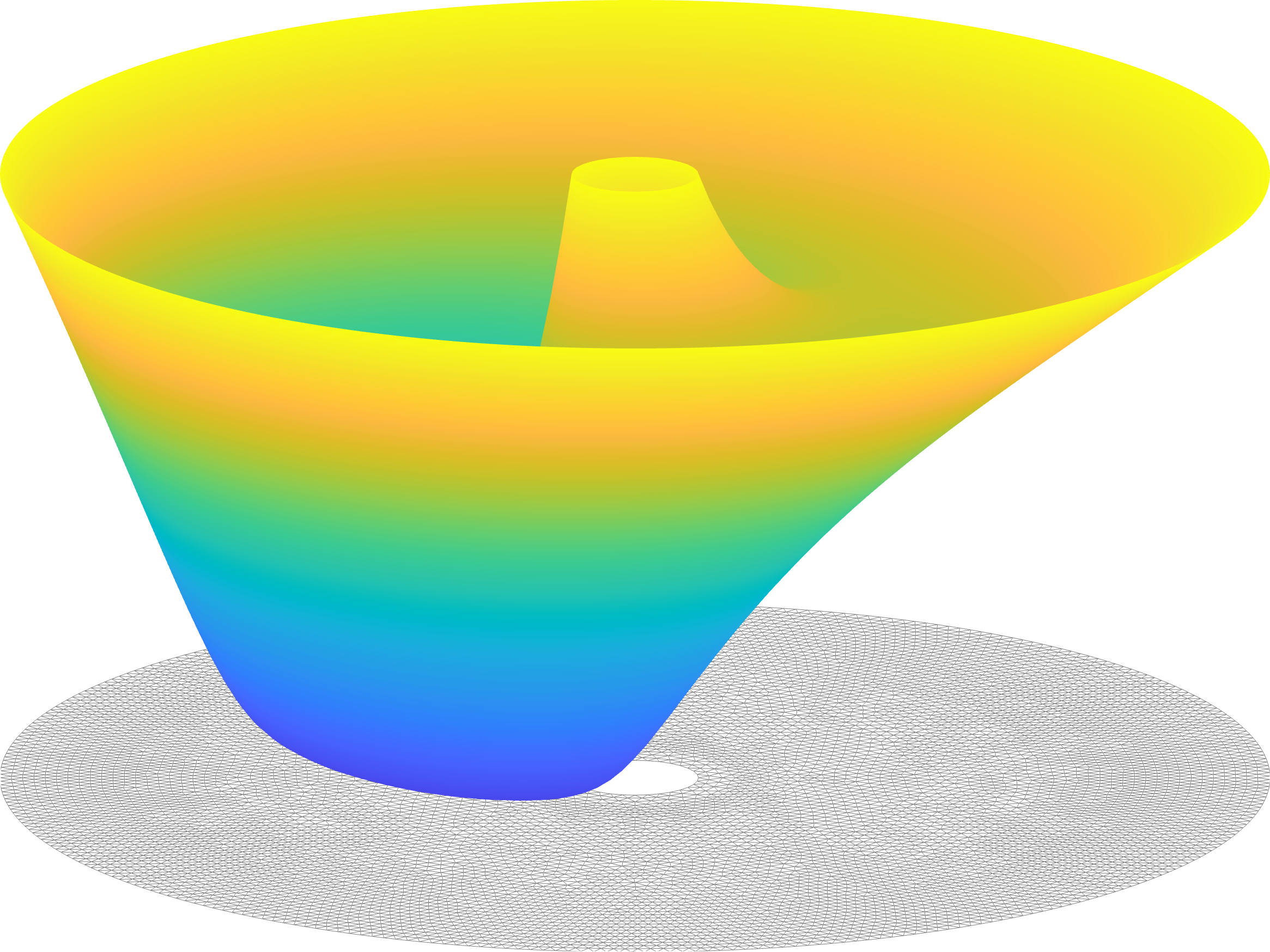};

\nextgroupplot[axis equal image, axis lines=none, width=0.45\textwidth, xtick={0, 1}, ytick={0}, title={6}, xticklabels={}, yticklabels={}, xmin=0, xmax=1, ymin=-1, ymax=0]
\addplot []
graphics [xmin=0,xmax=1,ymin=-1,ymax=0] { 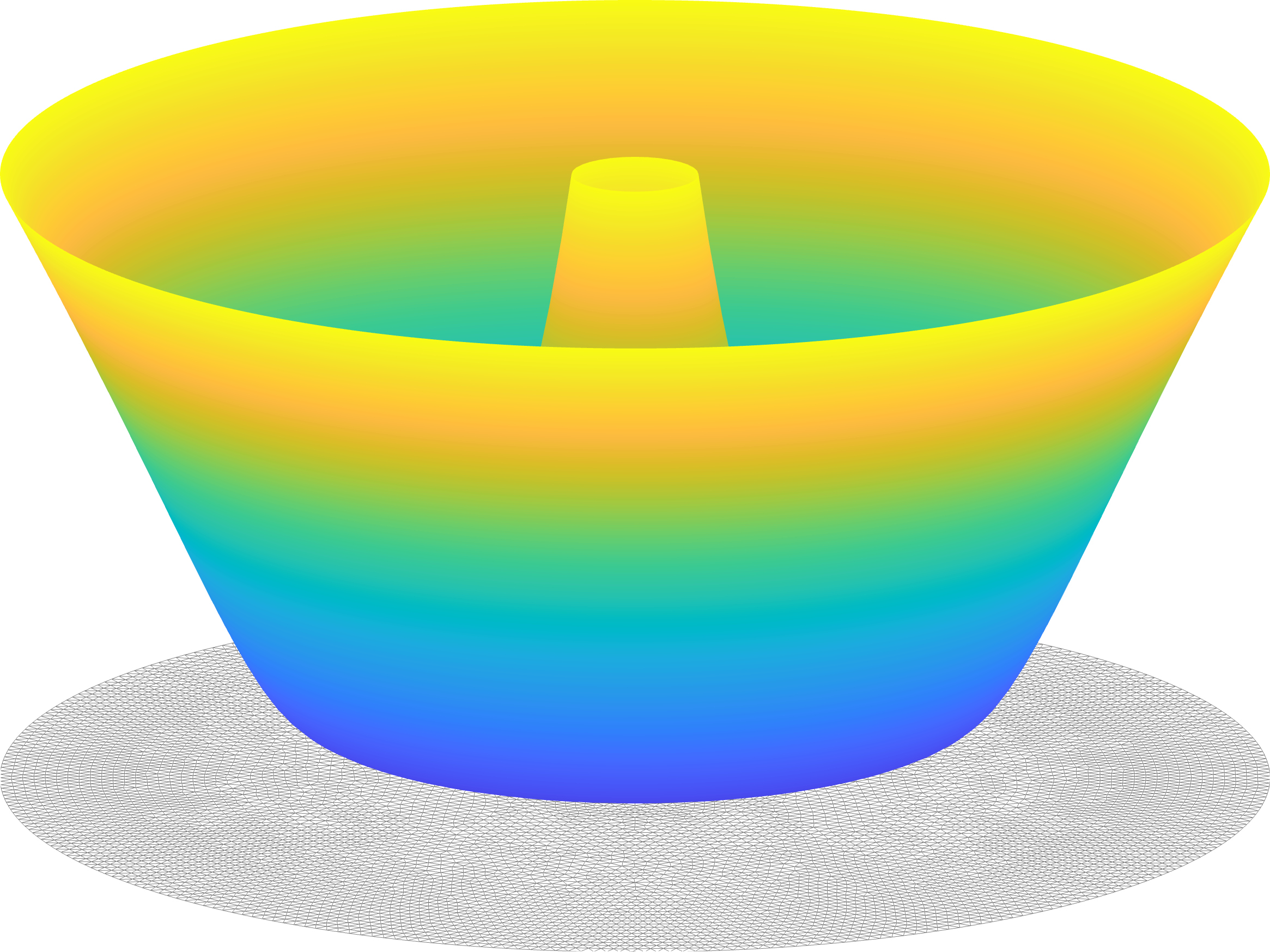};

\nextgroupplot[axis equal image, axis lines=none, width=0.45\textwidth, xtick={0, 1}, ytick={0}, title={7}, xticklabels={}, yticklabels={}, ylabel={$\varepsilon = 0.1$}, xmin=0, xmax=1, ymin=-1, ymax=0]
\addplot []
graphics [xmin=0,xmax=1,ymin=-1,ymax=0] { 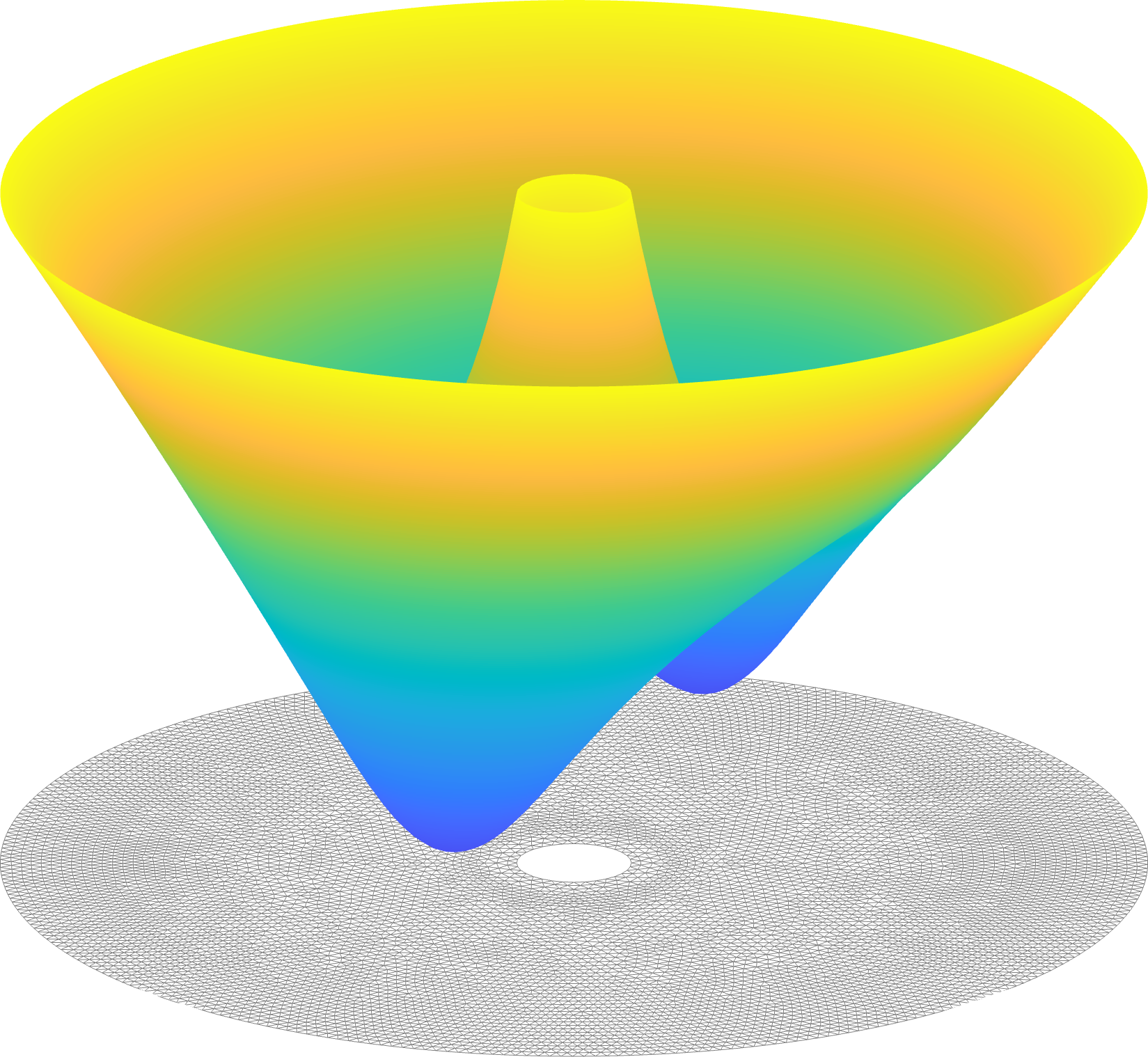};

\nextgroupplot[axis equal image, axis lines=none, width=0.45\textwidth, xtick={0, 1}, ytick={0}, title={8}, xticklabels={}, yticklabels={}, xmin=0, xmax=1, ymin=-1, ymax=0]
\addplot []
graphics [xmin=0,xmax=1,ymin=-1,ymax=0] { 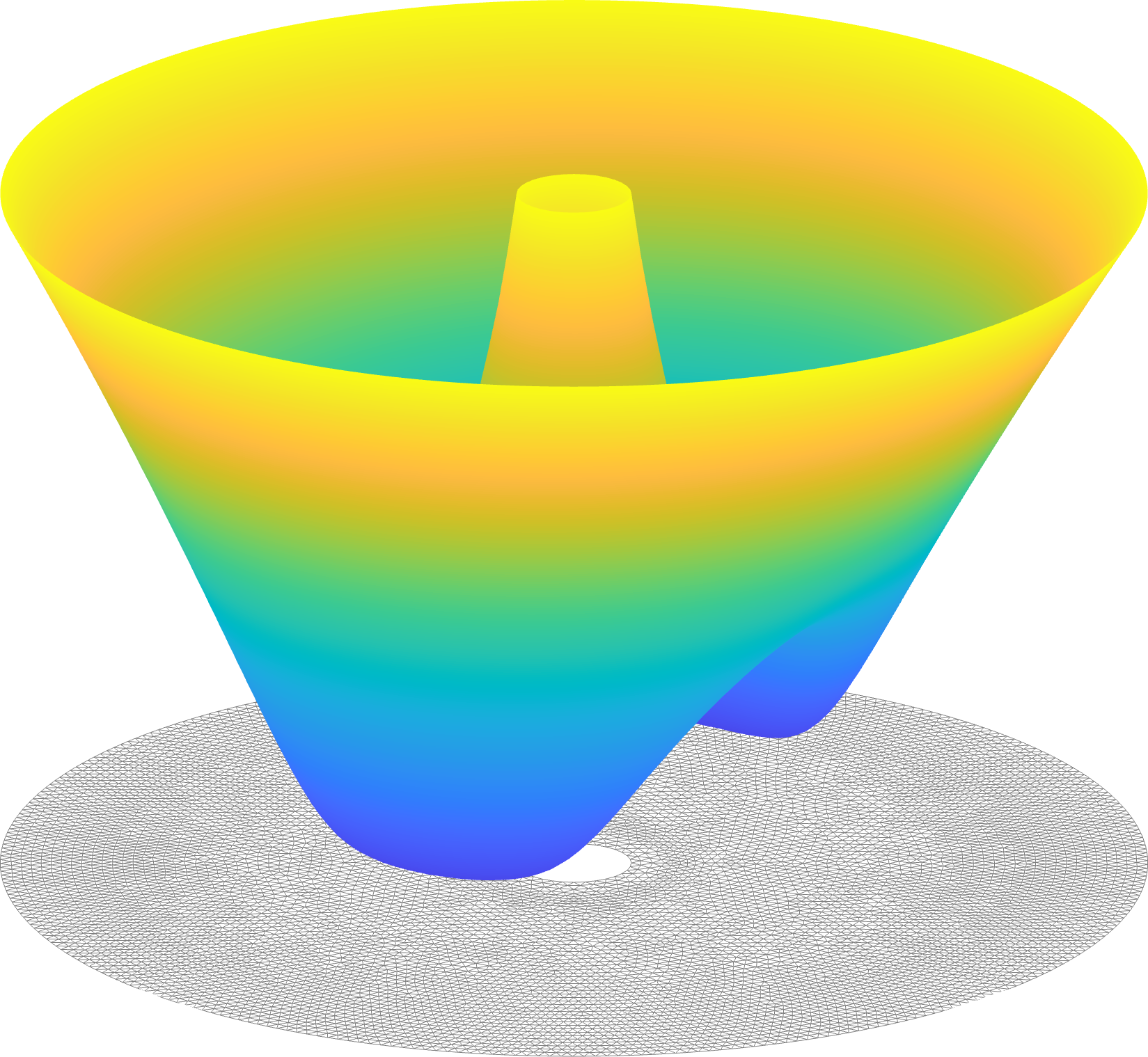};

\nextgroupplot[axis equal image, axis lines=none, width=0.45\textwidth, xtick={0, 1}, ytick={0}, title={9}, xticklabels={}, yticklabels={}, xmin=0, xmax=1, ymin=-1, ymax=0]
\addplot []
graphics [xmin=0,xmax=1,ymin=-1,ymax=0] { 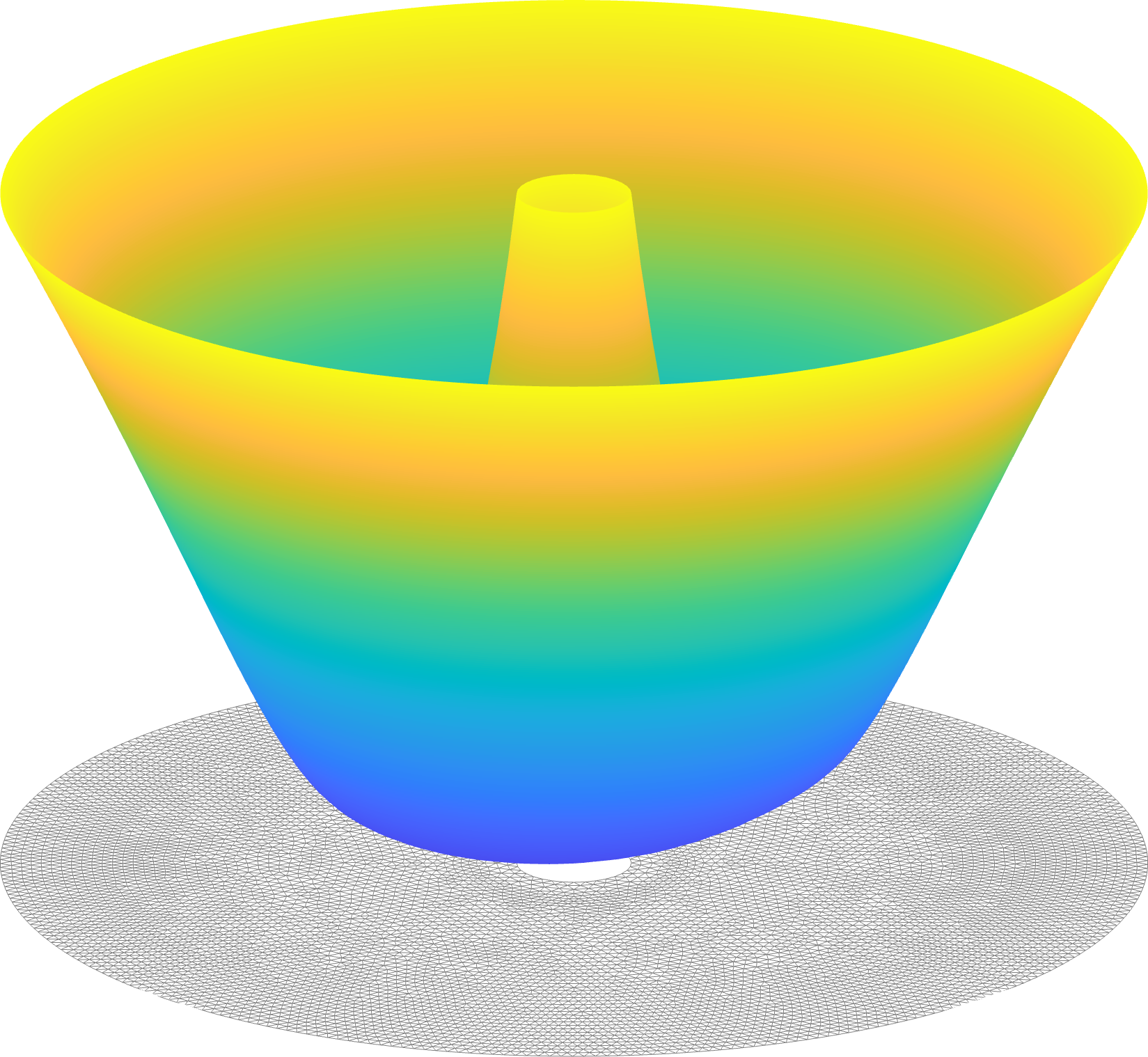};

\end{groupplot}\end{tikzpicture}
    }
\colorbarMatlabParula{-1}{-0.75}{-0.5}{-0.25}{0}
 
 \caption{ \revthree{ Solution of \eqref{eq:steady_problem} for $\eps=0.2$ and $\eps=0.1$ on the annular region $0.1<|\bx|<1$. Left: Bifurcation curves for $\|\bu\|_2$ (solid blue)(\ref{line:bif_curve_annulus_branch_1}) showing main branch and asymmetric branches. Right: Solutions and associated computational meshes for each of the labeled points for each branch. } }  \label{fig:asymBranches}

\end{figure}

\revthree{We now discuss the existence of these symmetry breaking bifurcations by analyzing the linearized system. Let us consider equation \eqref{mems_intro_eqns_a} for $\eps=0$ in the annulus $\Omega = \{ r(\cos\theta,\sin\theta) \in \mathbb{R}^2 \ | \  r\in(r_1,1), \ \theta \in (0,2\pi) \}$. 
Writing $u(r,\theta,t) = u^{\ast}(r) + e^{-\mu_k t}e^{ik\theta}\phi_k(r)$ for $|\phi_k| \ll1$, we find that
\begin{equation}
     \frac{\partial^2 \phi_k}{\partial r^2}  + \frac{1}{r}\frac{\partial \phi_k}{\partial r} + \left[\frac{2\lambda}{(1+u^{\ast})^3} - \frac{k^2}{r^2}\right] \phi_k = -\mu_k \phi_k,
\end{equation}
where $\phi_k(r_1) = \phi_k(1)=0$. If $\mu_k(\lambda) = 0$ for  $\phi_k(r)\neq0$ at some values of $k$ and $\lambda$, then a symmetry breaking bifurcation has been located. Sturm-Liouville theory guarantees a sequence $\{\mu_{k,1},\mu_{k,2},\mu_{k,3},\ldots\}$ with $\mu_{k,q} \to \infty$ as $q\to\infty$, though it remains to show that $\mu_{k,1}$ can be negative indicating an instability. We can define $V = \{ v\in L^2(\Omega) \ | \ v(\partial\Omega) = 0\}$ and write
\begin{equation}\label{eq:steadystate_radial1}
\mu_{k,1} = \inf_{\phi\in V} \frac{\int_{r_1}^1(\phi_r^2 - \frac{2\lambda}{(1+u^{\ast})^3}\phi^2 + \frac{k^2}{r^2}\phi^2) dr}{\int_{r_1}^1 \phi^2 r dr} \leq \frac{\int_{r_1}^1 (u^{\ast^2}_r - \frac{2\lambda}{(1+u^{\ast})^3}u^{\ast^2} + \frac{k^2}{r^2}u^{\ast^2}) r dr}{\int_{r_1}^1 u^{\ast^2} r dr} = \frac{\chi(u^{\ast})}{\int_{r_1}^1 u^{\ast^2} r dr}.
\end{equation}
On the upper branch we note that $\min_{\Omega}u^{\ast} \to -1$ as $\lambda\to0^{+}$, hence for any arbitrarily small $\alpha>0$, we have for some $p\in\mathbb{R}$ that
\begin{equation}\label{eq:unbounded}
\frac{2u^{\ast}}{1+u^{\ast}} \leq -p^2, \qquad u^{\ast}\leq -1 + \alpha.
\end{equation}
In other words, the term $2u^{\ast}/(1+u^{\ast})$ becomes arbitrarily large and negative as the branch is traversed. Using this fact, together with the identity $\int_{\Omega} |\nabla u^{\ast}|^2dx =-\lambda\int_{\Omega}u^{\ast}/(1+u^{\ast})^2 dx$, it was shown in \cite[Sec.~5]{Feng2005} that
\begin{equation}\label{eq:stability}
    \chi(u^{\ast}) = \int_{r_1}^1 (u^{\ast^2}_r - \frac{2 u^{\ast}}{(1+u^{\ast})}\frac{\lambda u^{\ast}}{(1+u^{\ast^2})} + \frac{k^2}{r^2}u^{\ast^2}) r dr \leq \Big(1-p^2 + \frac{k^2}{r_1^2 \nu_1}\Big) \int_{r_1}^{1}u^{\ast^2}_r rdr - M,
\end{equation}
where $M>0$ is a constant independent of $\lambda$ and $\nu_1$ is the principal eigenvalue of the Laplacian on $\Omega$ with Dirichlet boundary conditions. Hence, for any given mode $k$, we can see that for $\lambda$ sufficiently small ($p^2$ sufficiently large), that $\mu_{k,1}(\lambda)<0$ thus indicating a symmetry breaking bifurcation.

In the case $\eps>0$, we find that the bifurcation structure changes again and we now observe a finite number of symmetry breaking bifurcations. In Fig.~\ref{fig:asymBranches} we plot bifurcation diagrams and solution profiles for values $\eps=0.2$ and $\eps=0.1$. As this bifurcation parameter decreases from $\eps=0.2$ to $\eps = 0.1$, we observe that the number of symmetry breaking bifurcations increases from one to two. We hypothesise that for any positive $\eps>0$, there will be a finite number $k_c$ of symmetry breaking bifurcations, with $k_c\to\infty$ as $\eps\to0$.
}

\section{Discussion}\label{sec:discussion}

In this work we numerically characterized two dimensional solutions of the system \eqref{eq:steady_problem} over the range of parameters $\lambda > 0$, $\eps > 0$ and $m = 4$. In Section \ref{sec:fem} we outlined a finite element discretization that couples adaptive mesh refinement (Sec.~\ref{sec:adapt_href}) with pseudo arc-length continuation (Sec.~\ref{sec:newt_raph}) to trace out the family of solutions to \eqref{eq:steady_problem}. We provided analysis of post touchdown states over a range of values of $\eps$, including relatively large values $(\eps = 0.3)$ and relatively small values $(\eps = 0.01)$ on the unit disk domain (Sec.~\ref{sec:reg_unit_disk}) and the unit square (Sec.~\ref{sec:reg_unit_square}). Finally, in Section \ref{sec:symm_break} we explore asymmetric solutions in the annulus domain using the methodology outlined in Section \ref{sec:branch_switch}.

The adaptive mesh refinement algorithm provided dynamic relocation of points to areas of the domain with fine scale features, and coarsening to areas of the domain with constant or near constant solution. The pseudo arc-length continuation procedure allows for navigation of solution branches and for solution multiplicity to be identified in terms of the bifurcation parameter $\lambda$. The numerical method also showed flexibility towards non-trivial geometries, the annulus domain, and still maintained its robustness in adapting the mesh resolution to approximate different fine scale features that develop in different regions of the domain. Overall, adaptive mesh refinement coupled with pseudo arc-length continuation provides a robust method to explore post equilibrium states described by \eqref{eq:steady_problem}, particularly those featuring fine scale solution features. 

In terms of future work, deflation is one method which has the potential to improve upon the computational approach used here. The premise of deflation \cite{Farrell2015,Butten2022} is to factor out previously discovered solutions to the system so that only new solutions are obtained from future iterations. An advantage of this method over pseudo arc-length continuation is that it can readily compute isolated solutions. However, a challenge to effective implementation will be to incorporate adaptivity in the computation of new solutions which may have very different meshing requirements. A combination of moving mesh methods \cite{LD2017,Walsh2013,ansiotropic2015,DIPIETRO2018763,DiPietro2019} with deflation may be advantageous since degrees of freedom remain fixed through iterations while allowing for localized refinement.

\revone{The number of degrees of freedom increases significantly with the use of adaptive meshing. Efficiencies can be gained by implementing a mesh-free methodology to bypass forming the full Jacobian \cite{lust2000computation}. This approach can also be adapted to investigate the presence of periodic orbits \cite{sanchez2010multiple,kavousanakis2008timestepper}.}

The family of solutions to \eqref{eq:steady_problem} is rich and complex over all parameters and geometries. The numerical tool developed here can be extended to study problems related to \eqref{eq:steady_problem}, such as different geometries, other types of forcing functions, or other types of symmetry breaking configurations (topological changes). Potential future work could focus on multi-parameter continuation methods for both $\eps$ and $\lambda$ \cite{dankowicz2020multidimensional, henderson2002multiple, henderson2007higher}. An investigation of the full two parameter $(\eps,\lambda)$ bifurcation structure would start with a mature bifurcation software package (e.g. \textsc{MatCont} \cite{Matcont}) with efforts taken to minimize number of unknowns.

Finally, our numerical results also reveal that new analytical results are needed to characterize the solution multiplicity of \eqref{eq:steady_problem} and in particular to identify regions of $(\lambda,\eps)$ space where the system is bistable.



\section*{Acknowledgments} A.E.L acknowledges support under NSF DMS-1516753.

\bibliography{MEMSbib.bib}{}
\bibliographystyle{plain}


\end{document}